\documentclass[11pt]{article}
\usepackage{mathrsfs}
\usepackage{amssymb,amsmath,amsfonts,amsthm}
\usepackage{amsfonts}
\usepackage{graphicx,subfigure,epstopdf}
\usepackage[usenames]{color}
\usepackage{url}
\usepackage{algorithm,algorithmic}
\usepackage{dsfont,verbatim}
\usepackage[hidelinks,linktocpage]{hyperref}
\usepackage[margin=1.1in]{geometry}
\usepackage{enumerate,enumitem}
\usepackage[normalem]{ulem}
\usepackage{cite}
\usepackage{bm}
\usepackage{upgreek}
\allowdisplaybreaks

\newcommand{\nn}{\nonumber}

\newtheoremstyle{thmm}{1.5ex plus 1ex minus .2ex}{1.5ex plus 1ex minus.2ex}{\rmfamily}{}{\bfseries}{}{1em}{} \theoremstyle{thmm}
\newtheorem{theorem}{Theorem}[section]
\newtheorem{lemma}{Lemma}[section]
\newtheorem{proposition}{Proposition}[section]

\newtheorem{remark}{Remark}[section]
\renewenvironment{proof}[1][Proof]{\noindent\textit{#1. }
}{\hfill$\square$}

\newcommand{\vertiii}[1]
{{\left\vert\kern-0.25ex\left
\vert\kern-0.25ex\left\vert #1
    \right\vert\kern-0.25ex\right
\vert\kern-0.25ex\right\vert}}

\def\d{{\rm d}}
\def\E{{\mathbb E}}
\def\R{\mathbb{R}}

\def\N{\mathbb{N}}

\def\nu{n}

\def\z{{\bf z}}

\def\d{{\mathrm d}}

\def\R{{\mathbb R}}

\def\O{{\mathcal O}}
\def\H{{\mathcal H}}

\def\o{\Omega}

\newcommand\blfootnote[1]{%
  \begingroup
  \renewcommand\thefootnote{}\footnote{#1}%
  \addtocounter{footnote}{-1}%
  \endgroup
}

\title{\Large\bf Maximal ${\bf L^p}$-${\bf L^q}$ regularity of finite element methods for stochastic {parabolic} equations in polygonal domains}
\author{\normalsize
	Genming Bai
	\thanks{Department of Mathematics and Statistics,
		Old Dominion University, Norfolk, VA 23529, USA;
		and Department of Mathematics,
		University of Michigan, Ann Arbor, MI 48109, USA.
		Email address: gbai@odu.edu.}
	\and
	Jianbo Cui
	\thanks{Department of Applied Mathematics,
		The Hong Kong Polytechnic University,
		Hung Hom, Kowloon, Hong Kong.
		Email addresses: jianbo.cui@polyu.edu.hk and buyang.li@polyu.edu.hk.}
	\and
	Buyang Li$^\dagger$
}

\date{}

\begin{document}

\maketitle

\vspace{-10pt}

\begin{abstract}
We prove that a spatial semidiscrete finite element method preserves the stochastic maximal regularity of the heat equation {on} bounded convex polygonal and polyhedral domains in dimensions two and three. 
{A main difficulty is to establish suitable Green-function estimates
for stochastic convolutions, whose paths generally lack classical
time derivatives.}
Besides, the failure of the Burkholder–Davis–Gundy (BDG) inequality at the endpoint $p=\infty$ creates additional difficulties. Our approach leverages  a hidden symmetry structure of the discrete Green's function, which enables the use of a key integration-by-parts type formula. Then the discrete stochastic maximal $L^p$-$L^p$ regularity with $p\in [2, \infty)$ is obtained from the discrete local energy estimates and interpolation inequalities. {To extend further the} discrete stochastic maximal $L^p$-$L^q$ regularity with the full range of integrability exponents,  i.e., $p\in (2,\infty), q\in [2,\infty)$ or $p=q=2$,  we also show that the discrete convolution kernel satisfies the stochastic H\"ormander's condition allowing the application of a stochastic version of the Calder\'on–Zygmund theory. 
\\ 

\end{abstract}

\blfootnote{This work is partially supported by the National Natural Science Foundation of China (projects no. 12525111, 12301526, 12231003 and 12522119) and the Hong Kong Research Grants Council (projects no. RFS2324-5S03, 15301025, 15301818).}

\vspace{-12pt}

\setlength\abovedisplayskip{5pt}
\setlength\belowdisplayskip{5pt}

\section{Introduction}
\setcounter{equation}{0}
Maximal regularity is an important property of the abstract parabolic equation
\begin{align}\label{eq:PDE0}
	\left\{
	\begin{aligned}
		&\partial_t u +A u = f , \\ 
		&u(0)= 0  .
	\end{aligned}
	\right.	
\end{align}
It roughly says the two dominant parts on the left hand side of the parabolic equation, i.e., $\partial_t u$ and $A u$, are separately as regular as the source $f$. To be more rigorous, we say the operator $A$ has the maximal $L^p$-regularity if the $X$-valued solution $u(t)$ satisfies the a priori estimate
\begin{align}
	\| \partial_t u \|_{L^p(0,T; X)}
	+
	\| A u \|_{L^p(0,T; X)}
	\leq
	C
	\| f \|_{L^p(0,T; X)} ,
\end{align}
where the constant $C$ does not depend on {the terminal time} $T$; see \cite{KW04,Lunardi12}.
{More generally, finite-interval maximal regularity may have a constant depending on $T$, including in nonautonomous problems; see \cite{Pruss01}}.
In this paper, we will stick to the case of $A = -\Delta$ equipped with Dirichlet boundary condition and {$X=L^q(\Omega)$}.

Maximal $L^p$-regularity is of great importance in particular when dealing with nonlinear equations, i.e., $f = f(u,\nabla u)$, as it {enables}  the maximal regularity for both $\partial_t u$ and $A u$.  {This, in turn, allows for effective control of the nonlinear term  $f = f(u,\nabla u)$}. As a result, maximal $L^p$-regularity serves as {a} powerful tool in studying the well-posedness and regularity of solutions of nonlinear parabolic PDEs; see \cite{Amann1995,ClementPruss1992,LiYang2015,KW04,Lunardi12,Pruss01}.  Similarly, at the discrete level, the discrete maximal $L^p$-regularity  {facilitates the treatment of nonlinearities} (see \cite{ALL17,KLL18_focm,KL23}), and {enables the derivation of optimal  $L^p$ error estimates without relying on  inverse inequalities} (cf. \cite{KLL16,Li19,KL23}). This stronger error control leads to sharper results compared with $L^2$-based methods.

The above important implications of maximal regularity hold as well in the stochastic setting, where we consider the mild solution of the following SPDE
\begin{align}\label{eq:PDE1}
	\left\{
	\begin{aligned}
		&d u + A u \d t = f \d W, \\ 
		&u(0)= 0,
	\end{aligned}
	\right.	
\end{align}
with $W$ being the standard real-valued Brownian motion.
However, it is more difficult to establish the stochastic maximal regularity results due to the lack of regularity which stems from the roughness of the white noise. {Brownian paths are almost surely H\"older continuous of every order below $1/2$ (see \cite{Le16}), and stochastic forcing generally precludes classical time derivatives of the solution.} The loss of temporal regularity prevents us from using multiplier theory in the time variable, rendering the classical $R$-boundedness approach (cf. \cite{KW04}) for maximal regularity invalid. In \cite{Van12}, the authors established the $L^p$-$L^q$ stochastic maximal regularity result via a combination of $R$-boundedness and $H^\infty$-functional calculus.
Similar results for $C^1$ domains were obtained
 in \cite{Kim04a,Kim04b,LV21} whose proofs do not rely on the $H^\infty$-calculus. In \cite{CKLL18,CKL19}, a weighted version of stochastic maximal regularity is developed on an angular domain in $\R^2$.

Let the upright letter $\O$ be the sample space (we reserve the notation $\Omega\subset \mathbb R^N$ for the underlying physical domain) and $(\O, \mathscr F, (\mathscr F_t)_{t\geq 0}, \mathbb P)$ be the filtered probability space satisfying the usual conditions, i.e. the filtration $ (\mathscr F_t)_{t\geq 0}$ is complete and right-continuous (see \cite{Le16}).
We denote by $W$ the real-valued Brownian motion on this filtered probability space.
Let $q \in [2, \infty)$, and suppose that the operator $A$ admits a bounded $H^\infty$-calculus of angle less than $\pi/2$ on $L^q(\Omega)$. Let $(E(t))_{t\geq0}$ denote the bounded analytic semigroup on $L^q(\Omega)$ generated by $-A$, i.e. $E(t) = e^{-tA}$. {As established in \cite{Van12},} for any $(\mathscr F_t)_{t\geq 0}$-adapted process $f \in L^p (\O\times\R_+ ;L^q (\Omega))$ the stochastic convolution process
\vspace{-1mm}
\begin{align}
	u(t)
	=
	\int_0^t E(t-s) f(s) \d W_s, \quad t>0 ,
\end{align}
is well-defined in $L^q(\Omega)$, and takes values in the fractional domain $D(A^{1/2})$ almost surely. %
 For all $2 < p < \infty,$ we have the following stochastic maximal $L^p$-$L^q$ regularity estimate  \begin{align}
	\E \| A^{1/2} u \|_{L^p(\R_+; L^q(\Omega))}^p
	\leq
	C_{1,p,q}
	\E \| f \|_{L^p(\R_+; L^q(\Omega))}^p ,
\end{align}
where the constant $C_{1,p,q}$ may depend on $p$, $q$ and $\Omega$, but is
independent of $T$.
For $q = 2$ this estimate also holds with $p = 2$.
In particular, if we consider the heat equation, i.e. $A=-\Delta$, then
\vspace{-2pt}
\begin{align}\label{eq:cont-SMR}
	\E \| (-\Delta)^{1/2} u \|_{L^p(\R_+; L^q(\Omega))}^p
	\leq
	C_{1,p,q}
	\E \| f \|_{L^p(\R_+; L^q(\Omega))}^p .
\end{align}

Recently, several works have studied the maximal inequalities and uniform strong error for temporal discretization of the stochastic heat equation. In \cite{Van21}, the authors prove the discrete maximal inequality of stochastic convolution for a large class of time discretization schemes
(including splitting, implicit Euler, Crank--Nicolson, and other rational schemes). In \cite{KV24}, the uniform strong convergence for contractive time discretization schemes {is} shown for semilinear stochastic evolution equations with irregular Lipschitz nonlinearities. {We also refer to \cite{BCH19,BG19,KLL18,GN97} and references therein for the uniform strong error estimates for time discretization schemes of stochastic nonlinear heat equations.}

In this paper, we prove that the stochastic maximal $L^p$-$L^q$ regularity for the heat equation on convex polygonal domains in {$\R^N (N=2,3)$} is preserved {}under a spatially semidiscrete finite element method.
 Consider the following discrete stochastic convolution of the discrete semigroup $E_h(t)=e^{t\Delta_h}$, where $h>0$ corresponds to the mesh size and $\Delta_h$ is the discrete Dirichlet Laplacian associated to the underlying finite element space:
\begin{align}
	u_h(t)
	=
	\int_0^t E_h(t-s) f_h(s) \d W_s, \quad t>0 .
\end{align}
Here $f_h(\cdot)=P_h f(\cdot)$ with $P_h$ being the $L^2$-projection onto the finite element space.
{Our main result is the following stochastic maximal $L^p$-$L^q$ regularity estimate, which is established in Theorem \ref{thm:disc_max_reg2}:}   \vspace{-1.5mm}
\begin{align}
	\E \| (-\Delta_h)^{1/2} u_h \|_{L^p(0,T; L^q(\Omega))}^p
	\leq
	C_{2,p,q} 
	\E \| f_h \|_{L^p(0,T; L^q(\Omega))}^p ,
\end{align}
for all  $T>0$ and $p\in(2,\infty)$, $q\in[2,\infty)$ or $p=q=2$. 
{Here the constant $C_{2,p,q}$ is independent of $h$ and $T$,
but may depend on $p$, $q$, $\Omega$, the fixed polynomial degree,
and the shape-regularity and quasi-uniformity constants
of the mesh.}

In \cite{Kem18}, the authors show that $\Delta_h$ admits a bounded $H^\infty$-calculus if the lumped mass bilinear form is used (ensuring positivity) and the finite element triangulation satisfies the acuteness condition (ensuring discrete {maximum} principle, thus the contractivity in $L^\infty$).
{A more recent result \cite{LZ2026} shows that these finite element constraints can be removed given that $\Omega$ is convex and $\partial\Omega$ is of class $C^{2,0}$.} For more general situations, our strategy is to utilize Green's function methods, as introduced and developed in the deterministic cases (cf. \cite{Schatz80,Schatz95,Schatz98,Dem09,Ley09,Ley12,Li19,Li22,Li22b}).
Essentially, our proof only depends on the Gaussian behavior for Green's functions and does not rely on the $H^\infty$-calculus which is a stronger property than the desired stochastic $L^p$-$L^q$ regularity \cite[Remark 8.5(ii)]{LV21}.  {We also note that our proof can be extended to nonconvex domains with $C^{2,0}$ boundary in arbitrary dimensions. }

To address the difficulties in the stochastic setting, we developed the following techniques, which constitute the main novelties of this article:
\vspace{-1mm}
\begin{itemize}
	\item Hidden symmetry property of discrete Green's function:
	Based on the representation formula using Green's function, we have
	\begin{align}
		\nabla_{x_0} u_h(t,x_0) = \int_0^t \int_\Omega \nabla_{x_0} \Gamma_h(t-s,x,x_0) f(s,x) \d x \d W_s , \notag
	\end{align}
	where $\Gamma_h$ is the discrete Green's function of the discrete semigroup $e^{t\Delta_h}$.  Since we aim to utilize the local energy estimates, it is desirable to {transfer} the gradient with respect to the reference point $x_0$ to the other variable $x$ via {an} integration-by-parts argument.  Denote the Green's function of the semigroup $e^{-tA}$ by $G(t,\cdot,\cdot).$  Although $\Gamma_h(t,\cdot,\cdot)$ and $G(t,\cdot,\cdot)$ are symmetric, these symmetries generally do not extend to their derivatives, i.e., {in general,} $\nabla_{x_0} \Gamma_h(t,x,x_0) \neq \pm \nabla_x \Gamma_h(t,x,x_0)$ and $\nabla_{x_0} G(t,x,x_0) \neq \pm \nabla_x G(t,x,x_0)$. {In} some special cases, {such as when} $\Omega=\R^N$, the specific form of Green's function does enjoy translation invariance, i.e., $G(t,x,x_0) = G(t,x-x_0)$, and therefore $\nabla_{x_0} G(t,x,x_0) = - \nabla_x G(t,x,x_0)$. But this kind of integration-by-parts formula {is} not {available} for general domains. {Our key observation, stated in Lemma \ref{lemma:IbP}, provides a remedy:} 	\begin{align}
		- P_h (\nabla_x \tilde \delta_{x_0}(x))
		=
		P_h (\nabla_{x_0} \tilde \delta_{x_0}(x))
		=
		\nabla_{x_0} P_h (\tilde \delta_{x_0}(x)), \notag
	\end{align}
	where $\tilde \delta_{x_0}(\cdot)$ is a regularized delta function \cite[Lemma 2.2]{Tho00}.
	We have discovered  that, although the derivatives of $\Gamma_h$ need not be symmetric, the initial value $\Gamma_h(0)$ satisfies an integration-by-parts type formula after {applying} the $L^2$-projection $P_h$ at the front. This formula will essentially allow us to remove the derivative on the reference point $x_0$, thereby enabling the application of local energy estimates.
    \vspace{-1.5mm}
	\item Boundedness of $L^\infty_{x_0} L^2_t L^1_x$-norm of $F(t,x,x_0)$ (Lemma \ref{LemGm2}):  Let $F(t,\cdot,\cdot)$ denote the error of the Green's function in the gradient sense (see \eqref{err-green} for the definition).
	According to the Burkholder–Davis–Gundy (BDG) inequality, it is necessary to consider the $L^2$ norm of the temporal variable. This basically results from the fact that the It\^o integral is constructed by Hilbertian methods. As a consequence, in the local energy estimates, we establish the boundedness of $\| F(t,x,x_0) \|_{L^\infty_{x_0} L^2_t L^1_x}$, which is stronger than $L^\infty_{x_0} L^1_t L^1_x$-norm obtained in the deterministic setting \cite{Schatz98,Li19}.
    \vspace{-1.5mm}
	\item Smoothing property of the discrete semigroup $e^{t\Delta_h}$:
	In Step IV of the proof of Lemma \ref{LemGm2}, we develop the following discrete version of the smoothing property:
	For large enough $n$ independent of $h$,
	\begin{align}
		\| (-\Delta_h)^{-n} P_h \nabla_x \tilde\delta_{x_0} \|_{L^2}
		\leq
		C \| \tilde\delta_{x_0} \|_{L^1} . \notag
	\end{align}
	This discrete smoothing effect can be obtained by a careful duality argument and the stability of $P_h$.
   \vspace{-1.5mm}
	\item To extend the results from $L^p$-$L^p$ to $L^p$-$L^q$, one popular approach is based on the singular integral theory. 
    However, in the stochastic setting, {the trajectory does not admit classical time derivatives.} 
Consequently, methods that rely on time derivatives, including vector-valued singular integral theory in the time direction, are not applicable.	
	A recently developed stochastic version of singular integral theory in \cite[Theorem 8.2, Remark 8.3 (iv)]{LV21} {provides a way forward}. {It suffices to show that the deterministic kernel satisfies a stochastic H\"ormander's condition; see Proposition \ref{prop:stochastic-extrapolation}. This can be checked by the discrete smoothing estimate in
	Lemma \ref{lem:pt_est2}.}
\end{itemize}
   \vspace{-1mm}
This article is organized as follows.
In Section \ref{sec-2}, we state the main theorem of this paper. In Section \ref{sec-3}, we {present} the preliminary results of Green's function estimates and local energy estimates. Finally, Section \ref{sec-4} proves the discrete stochastic maximal regularity.

\section{The main results}
\label{sec-2}
\setcounter{equation}{0}

Let $(\O, \mathscr F, (\mathscr F_t)_{t\geq 0}, \mathbb P)$ be the underlying filtered probability space satisfying the usual conditions and $W$ be the real-valued Brownian motion on this filtered space. The Burkholder--Davis--Gundy (BDG) inequality of the following form (cf. \cite[Proposition 2.1]{Van12}) will be used later.
\begin{lemma}[The BDG inequality]\label{lemma:BDG}
	Let $p \in (1,\infty)$ and $q \in (1,\infty)$ be fixed. 
    {For every strongly measurable $(\mathscr F_t)_{t\geq 0}$-predictable process} $f: \O\times\R_+ \rightarrow L^q(\Omega)$ we have
	\begin{align}
		 c_{p,q}^p \E \| f \|_{L^q(\Omega; L^2(\R_+))}^p
		\leq
		\E \bigg\| \int_0^\infty f \d W_s \bigg\|_{L^q(\Omega)}^p
		\leq
		C_{p,q}^p \E \| f \|_{L^q(\Omega; L^2(\R_+))}^p . \notag
	\end{align}
	Here $ c_{p,q}$ and $C_{p,q}$ are positive constants.
\end{lemma}

On a bounded domain $\Omega$ of arbitrary dimensions $N$, the Dirichlet Laplacian $\Delta$ on $\Omega$ has a bounded $H^\infty$-calculus of angle less than $\pi/2$ (see \cite[Example 2.8]{Van12}). We recall the following continuous stochastic maximal regularity result obtained in \cite[Theorem 1.1]{Van12}.
\begin{proposition}[Continuous maximal regularity of stochastic convolution]\label{thm:max_reg}
	{\it
Let $p\in(2,\infty)$, $q\in[2,\infty)$ or $p=q=2$.
        {For every strongly measurable $(\mathscr F_t)_{t\geq 0}$-predictable process} $f \in L^p (\O\times\R_+ ;L^q (\Omega))$ the stochastic convolution process
		\begin{align}\label{eq:stoh_conv_cont}
			u(t)
			=
			\int_0^t e^{(t-s)\Delta} f(s) \d W_s, \quad t>0 ,
		\end{align}
		is well defined in $L^q(\Omega)$, takes values in the fractional domain $D((-\Delta)^{1/2})$ almost surely and we have the stochastic maximal $L^p$-$L^q$ regularity estimate
		\begin{align}\label{eq:CSMR}
			\E \| (-\Delta)^{1/2} u \|_{L^p(\R_+; L^q(\Omega))}^p
			\leq
			C_{1,p,q}
			\E \| f \|_{L^p(\R_+; L^q(\Omega))}^p.
		\end{align}
	}
\end{proposition}
{Next, we record a semigroup extrapolation result via the stochastic singular integral theory (cf. \cite[Definition 8.1, Theorem 8.2, Remark 8.3 (iv)]{LV21} with the Muckenhoupt weight function $w$ set to constant $1$ and $H=\R$ therein):
	\begin{proposition}[\!\!{\cite[Theorem~8.2 and Remark~8.3(iv)]{LV21}}]
		\label{prop:stochastic-extrapolation}
		Let $X,Y$ be UMD spaces of type $2$, with $Y\hookrightarrow X$,
		and let $A$ be sectorial of angle less than $\pi/2$ on $X$.
		Suppose that the following smoothing estimate holds:
		\[
		\|e^{-tA}\|_{\mathcal L(X,Y)}
		\leq Mt^{-1/2}, \qquad t>0.
		\]
		For a real-valued Brownian motion, set
		$S_AG(t)=\int_0^t e^{-(t-s)A}G(s)\,\d W_s$.
		If $S_A$ is bounded from
		$L^{p_0}(\O\times\R_+;\gamma(\R,X))$ into
		$L^{p_0}(\O\times\R_+;Y)$ for some $p_0\in[2,\infty)$,
		then the same boundedness holds for every $p\in(2,\infty)$.
	\end{proposition}
}
We denote by $\Delta_h$ the discrete Dirichlet Laplacian associated to the finite element space $S_h^0(\Omega):=S_h(\Omega)\cap H_0^1(\Omega)$, where $S_h(\Omega)$ is the $H^1$-conforming finite element space consisting of piecewise polynomials of degree $r\geq 1$ subject to a quasi-uniform triangulation of  a convex polygon or polyhedron $\Omega$ in $\mathbb R^N, N=2,3$, with mesh size $0<h<h_0$ for some constant $h_0$. Our main result is stated below. {Let $\mathcal T_h$ denote this triangulation. Throughout this paper, the mesh size is defined by the maximal diameter of all elements \[ h:=\max_{\tau\in\mathcal T_h}\operatorname{diam}(\tau).\]}
\vspace{-4mm}
\begin{theorem}[Discrete stochastic maximal regularity of semidiscrete convolution via FEM]\label{thm:disc_max_reg2}
	{\it Let $T>0,  p\in(2,\infty)$, $q\in[2,\infty)$ or $p=q=2$.
        {For every strongly measurable $(\mathscr F_t)_{t\geq 0}$-predictable process} $f \in L^p (\O\times\R_+ ;L^q (\Omega))$ the discrete stochastic convolution process
		\begin{align}\label{eq:stoh_conv_disc}
			u_h(t)
			=
			\int_0^t e^{(t-s)\Delta_h} f_h(s) \d W_s \quad t>0 ,
		\end{align}
		satisfies the following estimate
		\begin{align}
			\E \| {{(-\Delta_h)^{1/2}}} u_h \|_{L^p(\R_+; L^q(\Omega))}^p
			\leq
		C_{2,p,q} 
			\E \| f_h \|_{L^p(\R_+; L^q(\Omega))}^p ,
		\end{align}
	where $f_h=P_h f$ is the spatial $L^2$ projection of $f$ onto the finite element space $S_h^0(\Omega)$ and the constant $C_{2,p,q}$ is independent of $h$. Consequently, restriction to $(0,T)$ gives the same estimate with a constant independent of $T>0$.
	}
\end{theorem}

{\begin{remark}\upshape
		We note that stochastic maximal regularity for temporal discretizations was recently studied in \cite{EV2026}. Combining our spatially semidiscrete result with \cite[Theorem~1.2]{EV2026} would yield the corresponding stochastic maximal regularity for fully discrete finite element schemes. We do not carry out this  additional verification here.
	\end{remark}

    \begin{remark}\upshape
		For simplicity, the main theorem is formulated for a real-valued
		Brownian motion. It also extends to an $\H$-cylindrical Brownian
		motion $W_\H$, where $\H$ is a separable Hilbert space.
		
		Indeed, set
		\[
		X_h=(S_h^0(\Omega),\|\cdot\|_{L^q}),\qquad
		A_h=-\Delta_h,\qquad E_h(t)=e^{-tA_h}.
		\]
		Applying \cite[Propositions~3.4 and~3.5]{LV21} to the kernel
		\[
		K_h(t,s)=\mathds{1}_{\{0<s<t\}}A_h^{1/2}E_h(t-s)
		\]
		and the scalar-noise estimate in
		Theorem~\ref{thm:disc_max_reg2}, we obtain
		\[
		\|A_h^{1/2}u_h\|_{L^p(\O\times(0,T);L^q(\Omega))}
		\le C_{p,q}
		\|f_h\|_{L^p(\O\times(0,T);\gamma(\H,X_h))},
		\]
		where
		\[
		u_h(t)=\int_0^t E_h(t-s)f_h(s)\,\d W_\H(s),
		\]
		for every strongly measurable $(\mathscr F_t)_{t\geq 0}$-predictable integrand $f_h$, with the same range of exponents
		as in Theorem~\ref{thm:disc_max_reg2}.
		
		The constant $C_{p,q}$ is independent of $h$, $T$, and the noise space
		$\H$, since the UMD and type $2$ constants of $X_h$ are
		uniformly controlled by those of $L^q(\Omega)$.
		Here $f_h=P_h\circ f$, and
		$\gamma(\H,L^q(\Omega))\cong L^q(\Omega;\H)$;
		see \cite[Remark~2.2]{Van12}.
	\end{remark}

}

\section{Preliminaries}
\label{sec-3}
\setcounter{equation}{0}
In this section, we introduce the underlying function spaces,  estimates related to Green's functions and the dyadic decomposition. {Notably,}  in Section \ref{sec:Grenn}, we discover a hidden symmetry of the discrete Green's function $\Gamma_h$ which allows us to perform an integration-by-parts type argument (see Lemma \ref{lemma:IbP}). We also develop the estimates of Green's functions associated to delta functions with spatial derivatives.

\subsection{Function spaces}
\label{sec:func}

We use the conventional notations of Sobolev–Slobodeckij spaces $W^{s,q}(\Omega)$, $s\ge 0$ and $1\leq q\leq\infty$, 
with the abbreviations $L^q=W^{0,q}(\Omega)$, $W^{s,q}=W^{s,q}(\Omega)$ and $H^s:=W^{s,2}(\Omega)$. 
The notation $H^{-s}(\Omega)$ denotes the dual space of $H^s_0(\Omega)$ (the closure of $C^\infty_0(\Omega)$ in $H^s(\Omega)$). We denote  $W^{-s,q}:=(W_0^{s,q'})^\prime, q'=\frac{q}{q-1}$, for $s>0$ and $q\in (1,\infty)$.
{For $\eta\in(0,1)$, we define the spatial H\"older seminorm
	and norm by
	\[
	|g|_{C^{1,\eta}(\Omega)}
	:=
	\sup_{\substack{z,\bar z\in\Omega\\z\ne\bar z}}
	\frac{|\nabla g(z)-\nabla g(\bar z)|}{|z-\bar z|^\eta},
	\qquad
	\|g\|_{C^{1,\eta}(\Omega)}
	:=
	\|g\|_{W^{1,\infty}(\Omega)}
	+
	|g|_{C^{1,\eta}(\Omega)}.
	\]}
\vspace{-1mm}
{Let $T>0$ and $\mathcal Q=(0,T)\times\Omega$.
	For $p,q\in[1,\infty]$ and $s\geq0$, we use the notation
	\begin{align}\label{eq:Bochner}
		\|g\|_{L^pW^{s,q}(\mathcal Q)}
		:=
		\|g\|_{L^p(0,T;W^{s,q}(\Omega))}
		=
		\big\|\|g(\cdot)\|_{W^{s,q}(\Omega)}\big\|_{L^p(0,T)}
	\end{align}
	for strongly measurable functions
	$g:(0,T)\to W^{s,q}(\Omega)$.
	More generally, mixed norms are defined using nested Bochner
	spaces (cf. \cite[Chapters 1 and 2]{HVVW2016book}), with subscripts indicating the corresponding variables.
	For example, given $g(x,y,z)$ defined on the product space $\Omega_x\times\Omega_y\times\Omega_z$, $p,q\in[1,\infty]$ and $r\in(0,1)$, we write
	\[
	\|g\|_{L^p_xW^{1,q}_yC^{1,r}_z(\Omega_x\times\Omega_y\times\Omega_z)}
	:=
	\|g\|_{L^p(\Omega_x;\,W^{1,q}(\Omega_y;\,C^{1,r}(\Omega_z)))}.
	\]
	The spaces are nested from right to left, and each $L^\infty$
	norm is understood as an essential supremum.
	We omit the variable subscripts and domain arguments whenever
	they are clear from the context as in \eqref{eq:Bochner}.
	
	On subdomains, we use restriction norms.
	For $D\subset\Omega$, $1\leq q\leq\infty$, and $s\geq0$, we set
	\begin{equation}\label{Def-HsD}
		\|g\|_{W^{s,q}(D)}
		:=
		\inf_{\substack{\tilde g\in W^{s,q}(\Omega)\\
				\tilde g|_D=g}}
		\|\tilde g\|_{W^{s,q}(\Omega)}.
	\end{equation}
	Similarly, for a space-time subdomain $Q\subset\mathcal Q$, we define
    \begin{equation}\label{DefLpX}
		\|g\|_{L^pW^{s,q}(Q)}
		:=
		\inf_{\substack{\tilde g\in L^p(0,T;W^{s,q}(\Omega))\\
				\tilde g|_Q=g}}
		\|\tilde g\|_{L^p(0,T;W^{s,q}(\Omega))}.
	\end{equation}
	In each case, the infimum is taken over all extensions of $g$
	in the indicated ambient space.
	The same restriction convention also applies to
	the H\"older norms. 
}

We use the abbreviations 
\begin{align}\label{inner-products}
	(\phi,\varphi):=\int_\Omega \phi(x)\varphi(x)\d x,\qquad
	[u,v]:=\int_0^T\int_{\Omega} u(t,x)v(t,x)\d x \d t ,
\end{align} 
and denote $w(t)=w(t,\cdot)$ 
for any function $w$ defined on ${\mathcal Q}$. The notation $1_{0<t<T}$ will denote the characteristic function of the time interval $(0,T)$, i.e., $1_{0<t<T}(t)=1$ if $t\in(0,T)$ while $1_{0<t<T}(t)=0$ if $t\notin(0,T)$.

For convex polygonal domains in $\R^N, N=2,3$, we use {a parameter} $\alpha>1$ to {characterize the smoothness of domains} via the following lemma (cf. \cite[Corollary 18.18]{Dauge88}).

\begin{lemma}\label{lemma:reg_domain}
	Let $\Omega$ be a convex polygon or polyhedron in $\R^N, N=2,3$. Then there exists a positive constant $\alpha > 1$ (depending on the domain $\Omega$) such that the solution of the Poisson equation
      \vspace{-1mm}
	\begin{align*}
	\begin{cases}
			\Delta \varphi  = f &\text{in $\Omega$} , \\
			\varphi = 0  &\text{on $\partial \Omega$}
			,
		\end{cases} \end{align*}
	satisfies
	$$
	\| \varphi \|_{H^{1+\alpha}}
	\leq
	C
	\| f \|_{H^{-1+\alpha}} .
	$$
\end{lemma}
\begin{remark}
	For general polygonal domains,  {Lemma \ref{lemma:reg_domain}} holds for some $\alpha>1/2$. We refer to \cite[Corollary 18.13]{Dauge88}. Also see the discussion after \cite[Lemma 4.2]{Li19}. 
\end{remark}

\subsection{Properties of the finite element space}
\label{sec:hypo}
Throughout this paper, we assume that the conditions of Theorem \ref{thm:disc_max_reg2} hold. 
{The Ritz projection $R_h:H_0^1(\Omega)\to S_h^0(\Omega)$
is defined by
\[
(\nabla R_hv,\nabla\chi_h)
=
(\nabla v,\nabla\chi_h)
\qquad
\forall\,v\in H_0^1(\Omega),\quad
\forall\,\chi_h\in S_h^0(\Omega).
\]
}
For any subdomain $D\subset\Omega$, we denote by $S_h(D)$ the space of  functions in $S_h(\Omega)$ restricted to the domain $D$ with  $h>0$, and denote by $S_h^0(D)$ the subspace of $S_h(D)$ consisting of functions which equal zero outside $D$. For any given subset $D\subset\Omega$, we denote $B_d (D) =\{x\in\Omega: \operatorname{dist}(x,D)\leq d\}$ for $d>0$. On a quasi-uniform triangulation of the domain $\Omega$, there exists a positive  constant $K>0$ such that the triangulation and the corresponding finite element space $S_h(\Omega)$ possess the following properties (cf. \cite[Appendix]{Schatz95}). Here $K$ is independent of the subset $D$ and $h$.

\begin{enumerate}[label={\bf (P\arabic*)},ref=\arabic*]\itemsep=5pt

\item {\bf Quasi-uniformity:}

\noindent For all triangles (or tetrahedra) $\tau_{par}^h$ in the partition,
the diameter $h_{par}$ of $\tau_{par}^h$ and the radius $\rho_{par}$
of its inscribed ball satisfy
$$
K^{-1}h\leq \rho_{par} \leq h_{par}\leq {h} ,
$$
where $par$ denotes the indices of the triangles (or tetrahedra) in the partition.

\item {\bf Inverse inequality:}

\noindent If $D$ is a union of elements in the partition, then
\begin{align*}
	\|\chi_h\|_{W^{l,p}(D)}
	\leq K h^{-(l-k)-(N/q-N/p)}\|\chi_h\|_{W^{k,q}(D)} ,
	\quad\forall\,\,\chi_h\in S_h(D),
\end{align*}
for $0\leq k\leq l\leq 1$ and 
$1\leq q\leq p\leq\infty$.  

\item {\bf Local approximation and super-approximation:} 

\noindent There exists an operator $I_h:H^1(\Omega)\rightarrow S_h( \Omega)$ with the following properties:

\begin{enumerate}[label=(\arabic*),ref=\arabic*]\itemsep=5pt
	\item 
	For $v\in W^{l,p}(\Omega ),l\in [1,2], p\in[1,\infty]$ such that $lp>N$, the following estimates hold
	:
	\begin{align*}
		&\|v-I_hv\|_{L^p} +h \|\nabla(v- I_hv)\|_{L^p} \leq Kh^l \|v\|_{W^{l,p}} 
		.
	\end{align*} 
	
	\item 
	If $d\geq 2h$ then the value of $I_h v$ in $D$ depends only on the value of $v$ in $B_d(D)$. 
	If $d\ge 2h$ and supp$(v)\subset \overline D$, then $I_hv\in S_h(B_{d}(D))$.

\item 
 Let  $g$ be a function on $\Omega$. 
If $d\geq 2h$, $g=0$ outside $D$
and $\sup\limits_{x\in D}|\partial^\beta g(x)|\leq Cd^{-|\beta|}$
for every multi-index $\beta$,
then  
\begin{align*}
	&\psi_h\in S_h(B_{d}(D))\implies I_h(g\psi_h)\in S_h(B_{d}(D)),\\
	&\|g\psi_h-I_h(g\psi_h)\|_{L^2}
	+h\|g\psi_h-I_h(g\psi_h)\|_{H^1}
	\leq K  h d^{-1}
	\|\psi_h\|_{L^2(B_{d}(D))} .
\end{align*}

\item  Let {$g$} be a function on $\Omega$. If $d\ge 2h$ and $g\equiv 1$ on $B_{d}(D)$, then $I_h(g\psi_h)=\psi_h$ on $D$. 

\item $I_hv\in S_h^0(\Omega)$ whenever
$v\in H_0^1(\Omega)$.

\end{enumerate}

\end{enumerate}

\subsection{Green's functions}
\label{sec:Grenn}
{The parabolic problems in this subsection
are subject to homogeneous Dirichlet boundary conditions in the
physical variable $x$.}
The heat equation on the domain $\Omega$ reads
\begin{align}\label{eq:HE}
	\left\{\begin{array}{ll}
		\partial_t u(t, x) - \Delta u(t, x)=f(t,x)
		&\mbox{in}\,\,\,(0,T)\times \Omega,\\
		u(0,x)=u_0(x)
		&\mbox{in}\,\,\,\Omega,
	\end{array}\right.
\end{align}
whose spatial semidiscrete finite element scheme is: Solving for $u_h\in S_h^0(\Omega)$ such that
\begin{align}\label{eq:HE_FEM}
	\left\{\begin{array}{ll}
		(\partial_t u_h(t,x),v_h(x))
		+ (\nabla u_h(t,x),\nabla v_h(x))
		= (f_h(t,x), v_h(x)) ,
		&\forall \, v_h\in S_h^0(\Omega), \\[5pt]
		u_h(0,x)=  u_{h,0},
	\end{array}\right.
\end{align}
where $f_h = P_h f$ and  $u_{h,0}=P_h u(0).$

We denote by $G(t,x,x_0)$ the Green's function
of equation \eqref{eq:HE} which satisfies
\begin{align}\label{GFdef}
	\left\{\begin{array}{ll}
		\partial_t G(t,x,x_0)-\Delta_x G(t,x,x_0)=0
		&\mbox{in}\,\,\, (0,T)\times \Omega,\\
		G(0,x,x_0)= \delta_{x_0}(x)
		&\mbox{in}\,\,\,\Omega .
	\end{array}\right.
\end{align}
The Green's function is symmetric, i.e. $G(t,x,y) = G(t,y,x)$. Moreover, it admits an {analytic} extension $G(z,x,y)$ on the open right half-plane $z\in\{t+is: t > 0{, s\in\R} \}$ (see \cite[Lemma 2]{Davies97}) satisfying the following pointwise estimate (cf. \cite[Theorem 3.4.8]{Davies89}).
\begin{lemma}\label{lemma:guass-est-c}
	{The analytic complex time heat kernel satisfies the following pointwise estimate.
	Fix $\theta_0\in(0,\pi/2)$. There is a constant $C_{\theta_0}>0$ such that    \begin{align}\label{eq:guass-est-c}
			|G(z,x,y)|\le C_{\theta_0}|z|^{-N/2}
			\exp\!\left(-\frac{|x-y|^2}{C_{\theta_0}|z|}\right)
		\end{align}
		for $z\ne0$ with $|\arg z|\le\theta_0$ and $x,y\in\Omega$.}
\end{lemma}

Then Cauchy's integral formula says that for all real time $t > 0$
\begin{equation}\label{Gaus-form}
	\partial_t^kG(t,x,y)
	=\frac{k!}{2\pi i}
	\int_{|z-t|=\frac{t}{2}}\,\, \frac{G(z,x,y)}{(z-t)^{k+1}}\d z ,
\end{equation}
which, together with \eqref{eq:guass-est-c}, yields the following Gaussian pointwise estimate for the time derivatives of Green's function:
\begin{align}
&|\partial_t^kG(t,x,x_0)|
	\leq 
	\frac{C_k}{t^{k+N/2}} e^{-\frac{|x-x_0|^2}{C_kt}}, &&
	\forall\, x,x_0\in \Omega,\,\,\, \forall\, t>0, \,\, k=0,1,2,\dots \label{GausEst1} .
\end{align} Here $C_k, k=0,1,2,\dots,$ are positive constants. 

As constructed in \cite[Lemma 2.2]{Tho00}, for any $x_0\in \bar\tau_{par}^h$ (where $\tau_{par}^h$ is an open {simplex} in the partition of $\o$), there exists a function $\tilde\delta_{x_0}\in C_0^\infty(\tau_{par}^h)$ such that
\begin{align*}
	\chi_h(x_0)=\int_{\Omega}\chi_h \tilde\delta_{x_0}\d x,
	\quad\forall\,\chi_h\in S_h(\Omega) ,
\end{align*}
and
\begin{align}
	&\|\tilde\delta_{x_0}\|_{W^{l,p}}
	\leq {C} h^{-l-N(1-1/p)}
	\quad\mbox{for}\,\,\,1\leq p\leq\infty,
	\,\,\, l=0,1,2,... , \label{reg-Delta-est}\\
	&\sup_{y\in \Omega} \int_\Omega |\tilde\delta_{y}(x)|\d x+
	\sup_{x\in \Omega}\int_\Omega |\tilde\delta_{y}(x)|\d y \le C .
	\label{reg-Delta-est-2}
\end{align}

Let $\delta_{x_0}$ denote the Dirac delta function centered at $x_0\in \Omega$. In other words, $\int_\Omega\delta_{x_0}(y)\varphi (y)\d y=\varphi(x_0)$ for arbitrary $\varphi\in C(\Omega)$.  
Then the discrete delta function is defined as
$$
\delta_{h,x_0} :=P_h  \delta_{x_0}=P_h \tilde\delta_{x_0} \qquad\forall \; x_0\in \Omega.
$$ 
It has been shown  (cf. \cite[Lemma 2.3]{Tho00} and \cite[Lemma 6.1]{Tho06}) that $\delta_{h,x_0}$ {decays} exponentially away from $x_0,$ i.e.,
\begin{align}\label{Detal-pointwise}
	& | \delta_{h,x_0}(x)|=|P_h \tilde\delta_{x_0}(x)|
	\leq {C} h^{-N} e^{-\frac{ | x - x_0 | }{C h }} ,
	\quad \forall\, x,x_0\in \Omega  ,
\end{align}
{and, similarly, using the locality and the scaling of $\nabla_x\tilde\delta_{x_0}(x)$, we have
\begin{align}\label{Detal-pointwise1}
	 |P_h \nabla_x \tilde\delta_{x_0}(x)|
	\leq {C} h^{-N-1} e^{-\frac{ | x - x_0 | }{C h }} ,
	\quad \forall\, x,x_0\in \Omega  .
\end{align}}
  \vspace{-0.5mm}
Recall the construction of the regularized delta function on the reference element $\hat \tau$: $\hat \delta_{x_0} (x) = \sum_{j=1}^{m_r}\hat\varphi_j(x_0) \hat\varphi_j(x) \hat w(x)$ where $\hat w\in C_0^{\infty}(\hat\tau)$ is a fixed and nonnegative function and $\hat\varphi_1, ..., \hat\varphi_{m_r}$ is an orthonormal basis for the local finite element space $S_h(\hat\tau)$ under the inner product $(\hat u, \hat v) = \int_{\hat\tau} \hat u \hat v \hat w \d\hat x$ (see \cite[Lemma 2.2]{Tho00}). 
From this formula, it can be seen that $\tilde \delta_{x_0} (x)$ and its derivatives are not (anti-)symmetric in $x$ and $x_0$. 
Our key observation is that, despite the lack of apparent (anti-)symmetry in $\nabla \tilde \delta$, we still have the following integration-by-parts type formula when the $L^2$ projection $P_h$ is applied at the front.
\begin{lemma}\label{lemma:IbP}For any $x\in \Omega$ and a.e. $x_0\in\Omega$,
	\begin{align}
		- P_h (\nabla_x \tilde \delta_{x_0}(x))
		=
		P_h (\nabla_{x_0} \tilde \delta_{x_0}(x))
		=
		\nabla_{x_0} (P_h \tilde \delta_{x_0}(x)) . \notag
	\end{align}
\end{lemma}
\begin{proof}
	For any given $x_0\in {\rm int}(\tau_{par}^h)$ where $\tau_{par}^h$ is an arbitrary element in the partition of $\Omega$, it suffices to prove 
	\begin{align}
		(-\nabla_x \tilde \delta_{x_0}(x), v_h)
		=
		(\nabla_{x_0} \tilde \delta_{x_0}(x), v_h) \quad\forall v_h\in S_h(\Omega) . \notag
	\end{align}
	Since by construction $\tilde\delta_{x_0}(x)\in C_0^\infty({\tau_{par}^h})$ and $\nabla_x v_h |_{\tau_{par}^h} \in  S_h({\tau_{par}^h})$, we use integration by parts to derive
	\begin{align}
		(-\nabla_x \tilde \delta_{x_0}(x), v_h)
		=
		(\tilde \delta_{x_0}(x), \nabla_x v_h)
		=
		\nabla v_h(x_0) . \notag
	\end{align}
	For the right-hand side,
	\begin{align}
	(\nabla_{x_0} \tilde \delta_{x_0}(x), v_h)
	=
		\nabla_{x_0} (\tilde \delta_{x_0}(x), v_h)
		=
		\nabla v_h(x_0). \notag	\end{align}\end{proof}

Let $\Gamma=\Gamma(t,x, x_0)$ 
be the regularized Green's function whose initial data is the regularized delta function $\tilde \delta_{x_0}$
\begin{align}\label{RGFdef}
	\left\{\begin{array}{ll}
		\partial_t\Gamma(t,x, x_0)-\Delta_x \Gamma(t,x, x_0)=0
		&\mbox{in}\,\,\,(0,T)\times \Omega,\\
		\Gamma(0,x ,x_0)=\tilde\delta_{x_0}(x)
		&\mbox{in}\,\,\,\Omega .
	\end{array}\right.
\end{align}
Then the regularized Green's function can be represented by 
\begin{align}\label{expr-Gamma}
	&\Gamma(t,x,x_0)=\int_\Omega  G(t,y,x)\tilde\delta_{x_0}(y)\d y=\int_\Omega  G(t,x,y)\tilde\delta_{x_0}(y)\d y  . 
\end{align}
From the representation \eqref{expr-Gamma} and \eqref{GausEst1}, one {can %
 derive} that the regularized Green's function $\Gamma$ also satisfies the Gaussian pointwise estimate in the region $\max(|x-x_0|, \sqrt t) \geq 2h$:
\begin{align}
	&|\partial_t^k\Gamma(t,x,x_0)|
	\leq 
	\frac{C_k}{t^{k+N/2}} e^{-\frac{|x-x_0|^2}{C_kt}}, \qquad
 \label{GausEstGamma0}
\end{align} 
with $k=0,1,2,\dots$.

Let $\Gamma_h=\Gamma_h(t,\cdot,  x_0)\in S_h^0(\Omega)$
be the discrete Green's function of the spatially semidiscrete scheme \eqref{eq:HE_FEM}
which approximates $\Gamma$,
defined by
\begin{align}\label{RGFdef1}
	\left\{\begin{array}{ll}
		\partial_t\Gamma_h(t,x, x_0)-\Delta_{h,x} \Gamma_h(t,x , x_0)=0
		&\mbox{in}\,\,\,(0,T)\times \Omega,\\
		\Gamma_h(0,x ,x_0)=\delta_{h,x_0}(x)
		&\mbox{in}\,\,\,\Omega ,
	\end{array}\right.
\end{align}
where the discrete Laplacian $\Delta_h g_h \in S_h^0(\Omega)$ %
{is defined as} the unique solution to
\begin{align}
	\int_\Omega \Delta_h g_h w_h {\d x}
	=
	- \int_\Omega \nabla g_h \nabla w_h {\d x}, \notag
\end{align}
for all $w_h\in S_h^0(\Omega)$.  Here, we use $\Delta_{h,x}$ to emphasize that the discrete Laplacian is acting on $x$, not $x_0$.
Then we differentiate the above evolution equation w.r.t. $x_0$ and get
\begin{align}\label{RGFdef2}
	\left\{\begin{array}{ll}
		\partial_t \nabla_{x_0} \Gamma_h(t,x, x_0)-\Delta_{h,x} \nabla_{x_0} \Gamma_h(t,x , x_0)=0
		&\mbox{in}\,\,\,(0,T)\times \Omega,\\
		\nabla_{x_0} \Gamma_h(0,x ,x_0)
		= - P_h (\nabla_{x} \tilde \delta_{x_0} )(x)
		&\mbox{in}\,\,\,\Omega ,
	\end{array}\right .
\end{align}
where we have used Lemma \ref{lemma:IbP} to derive
\begin{align}
	\nabla_{x_0}\delta_{h,x_0} 
	=
	\nabla_{x_0} (P_h \tilde\delta_{x_0} )
	= 
	- P_h (\nabla_{x} \tilde \delta_{x_0}) . \notag
\end{align}
{\begin{remark} 
	By the construction (cf. \cite[Lemma 2.2]{Tho00}),
		$\delta_{h,x_0}(x)$ is locally smooth in both variables
		$x$ and $x_0$ within the interiors of their respective mesh elements. 
		The time derivatives and elementwise spatial derivatives of $\Gamma_h$ can be interchanged there.
\end{remark}}
We define $H$ 
to be the regularized Green's function
associated to $\nabla_{x_0} \Gamma_h(\cdot,\cdot\, , x_0)$:
\begin{align}\label{RGFdef3}
	\left\{\begin{array}{ll}
		\partial_t H(t,x, x_0) - \Delta_{x} H(t,x , x_0)=0
		&\mbox{in}\,\,\,(0,T)\times \Omega,\\
		H(0,x ,x_0) = - \nabla_{x} \tilde \delta_{x_0}(x)
		&\mbox{in}\,\,\,\Omega .
	\end{array}\right.
\end{align}
This regularized Green's function $H$ can be represented by 
\begin{align}\label{expr-Gamma1}
	&H(t,x,x_0)=\int_\Omega  G(t,y,x) (- \nabla_{y} \tilde\delta_{x_0}(y))\d y  .
\end{align}
Analogous to \eqref{GausEstGamma0}, the following Gaussian pointwise estimate of $H$ holds in the region \\ $\max( {|x-x_0|}, \sqrt t) \geq 2h$ 
 using the representation \eqref{expr-Gamma1} and the estimates \eqref{reg-Delta-est} and \eqref{GausEst1}:
\begin{align}
	&|\partial_t^k H(t,x,x_0)|
	\leq 
	\frac{C_k}{h t^{k+N/2}} e^{-\frac{|x-x_0|^2}{C_kt}} \qquad
	  \label{GausEstGamma}
\end{align} 
with $k=0,1,2,\dots$.

\subsection{Dyadic decomposition of the domain ${\mathcal Q}=(0,{T})\times\Omega$}
\label{SecGF} 
In the proof of Theorem \ref{thm:disc_max_reg2}, {we} partition the domain ${\mathcal Q}$ into subdomains and {provide} estimates of the finite element solutions in each subdomain.  We present the dyadic decomposition of ${\mathcal Q}$ which   was introduced in \cite{Schatz98}. This subsection may be skipped by readers who are familiar with such dyadic decompositions.

For any integer $j$, we define $d_j=2^{-j}$.  For a given $x_0\in\Omega$, 
we let $J_1=1$, $J_0=0$ and $J_*$ be an integer satisfying $2^{-J_*}= C_*h$ 
with $C_*\geq \max \{16, 4K\}$ to be determined later, where $K$ is the constant of {quasi-uniformity} (see Section \ref{sec:hypo}).  {More precisely, we first choose a sufficiently large
constant $M\geq\max\{16,4K\}$, independently of $h$, and then
choose $C_*\in[M,2M)$ so that $J_*$ is an integer.
All subsequent estimates involving $C_*$ are understood
uniformly over this interval.} 
If 
\begin{align}\label{h-condition}
	h<1/(4C_*)  ,
\end{align}
then 
\begin{align}
	2\leq J_*=\log_2[1/(C_*h)]\leq \log_2(2+1/h).
\end{align} 
Let  
\begin{align*}
	&Q_*(x_0)=\{(t,x)\in \mathcal Q: \max ({|x-x_0|},t^{1/2})\leq d_{J_*}\}, 
	\\
	&\Omega_*(x_0)=\{x\in \Omega: |x-x_0| \leq d_{J_*}\} \, . 
\end{align*} 
We define 
\begin{align*} 
	&Q_j(x_0)=\{(t,x)\in {\mathcal Q}:d_j\leq
	\max (|x-x_0|,t^{1/2})\leq2d_j\}  &&\mbox{for}\,\,\, j\ge 1,
	\\
	&\Omega_j(x_0)=\{x\in \Omega: d_j\leq {|x-x_0|}  \leq2d_j\}   
	&&\mbox{for}\,\,\, j\ge 1,
\end{align*} 
and 
\vspace{-1mm}
\begin{align*} 
	&Q_0(x_0)= 
	{\mathcal Q}\big\backslash\big( \cup_{j=1}^{J_*}Q_{j}(x_0)\cup Q_*(x_0)\big) , 
	\\
	&\Omega_0(x_0) 
	=\Omega\big\backslash\big( \cup_{j=1}^{J_*}\Omega_{j}(x_0)\cup \Omega_*(x_0)\big).
\end{align*}
For $j<0$, we simply define
$Q_{j}(x_0)=\Omega_{j}(x_0)=\emptyset$.
For every integer $j\ge 0$, we define
\begin{align*}
	&\Omega_j'(x_0)=\Omega_{j-1}(x_0)\cup\Omega_{j}(x_0)\cup\Omega_{j+1}(x_0),
	\quad Q_j'(x_0)=Q_{j-1}(x_0)\cup Q_{j}(x_0)\cup Q_{j+1}(x_0), \\
	&\Omega_j''(x_0)=\Omega_{j-2}(x_0)\cup\Omega_{j}'(x_0)\cup\Omega_{j+2}(x_0),
	\quad
	Q_j''(x_0)=Q_{j-2}(x_0)\cup Q_{j}'(x_0)\cup Q_{j+2}(x_0).
\end{align*}
Then we have
\begin{align}\label{decomposition}
	&{\mathcal Q}=\bigcup^{J_*}_{j=0}Q_j(x_0)\,\cup Q_*(x_0)
	\quad\mbox{and}\quad
	\Omega=\bigcup^{J_*}_{j=0}\Omega_j(x_0)\,\cup \Omega_*(x_0) .
\end{align}{We also set $Q_*'=Q_*\cup Q_{J_*}$.
For each block index $i$, including $i=*$, let
$Q_i'(t):=\{x\in\Omega:(t,x)\in Q_i'\}$
denote its spatial section at time $t$.
}
We refer to $Q_*(x_0)$ as the ``innermost" set.
We shall write $\sum_{*,j}$ when the innermost set is included and
$\sum_j$ when it is not. When $x_0$ is fixed, if there is no ambiguity, 
we simply write
$Q_j=Q_j(x_0)$, $Q_j'=Q_j'(x_0)$, $Q_j''=Q_j''(x_0)$, 
$\Omega_j=\Omega_j(x_0)$, $\Omega_j'=\Omega_j'(x_0)$ and 
$\Omega_j''=\Omega_j''(x_0)$.

For $L^2$-based norms, we adopt the following abbreviations:
\begin{align}\label{Q-norm}
	&\|v\|_{k,D}=\biggl(\int_D\sum_{|\beta|\leq
		k}|\partial^\beta v|^2\d x\biggr)^{\frac{1}{2}},\qquad
	\vertiii{v}_{k,Q}=\biggl(\int_Q\sum_{|\beta|\leq
		k}|\partial^\beta v|^2\d x\d t\biggl)^{\frac{1}{2}} ,
\end{align}
for any subdomains $D\subset\Omega$ and $Q\subset (0,{T})\times\Omega$,  where $\beta$ is the multi-index for spatial variables.
We also simply write $\| v \|_{D}=\| v \|_{0,D}$ and $\vertiii{v}_{Q}=\vertiii{v}_{0,Q}$.
Throughout this paper, we denote by $C$ a generic positive constant that 
is independent of $h$, $x_0$ and $C_*$ (until $C_*$ is determined 
in Section \ref{ProofLm03}). To simplify the notation, we also denote $d_*=d_{J_*}$.

\subsection{Local energy estimates}
\label{sec:local}

In this subsection, we give two lemmas on the local energy estimates for Green's functions and the error equation. 
{\begin{lemma}\label{GFEst1} 
	{\it Given any $0\leq j\leq J_*$, the Green's function $G$ defined in \eqref{GFdef} satisfies the local high-order gradient estimate for $k\in \N$
		\begin{align}
			{\| \nabla \partial_t^k G(\cdot,\cdot,x_0) \|_{L^{\infty}L^\infty(Q_j(x_0))}
				\leq Cd_j^{-N-1-2k} }\label{GFest04},
		\end{align}
		and the local H\"older estimate
		\begin{align}
			{\|G(\cdot,\cdot,x_0) \|_{L^{\infty}C^{1,\eta}(Q_j(x_0))}
				\leq Cd_j^{-N-1-\eta} }\label{GFest03},
		\end{align}
		for some  $\eta\in (0,1)$ depending on the domain $\Omega$. Here and below the constant $C$ is independent of $x_0$.
		Consequently, the regularized Green's function 
		$H$ defined in \eqref{RGFdef3}
		satisfies the following estimates: 
		\begin{align}
			&d_j^{-4-\gamma+N/2}\|H(\cdot,\cdot,x_0)\|_{L^\infty(Q_j(x_0))}
			+d_j^{-4-\gamma}\|\nabla H(\cdot,\cdot,x_0)\|_{L^2(Q_j(x_0))}
			\nn\\
			&
			+d_j^{-4}\| H(\cdot,\cdot,x_0)\|_{L^2H^{1+\gamma}(Q_j(x_0))}  +d_j^{-2}\| \partial_{t} 
				H(\cdot,\cdot,x_0)\|_{L^2H^{1+\gamma}(Q_j(x_0))} \nn\\
			& 
			+\| \partial_{tt}H(\cdot,\cdot,x_0)\|_{L^2H^{1+\gamma}(Q_j(x_0))}
			\leq { C d_j^{{-N/2-5-\gamma}}} 			, \label{GFest01}
		\end{align}
		for all $0\leq\gamma \leq 1< \alpha$ (where $\alpha$ is determined in Lemma \ref{lemma:reg_domain}).
	}
	\end{lemma}
\begin{proof}
	The proof is divided into three steps. We first establish the local
	high-order gradient estimate \eqref{GFest04}. We then use this estimate,
	together with the spatial H\"older continuity of the elliptic Green's
	function, to prove \eqref{GFest03}. Finally, \eqref{GFest01} follows from
	the representation formula for $H$ and the argument used in
	\cite[Lemma~4.1]{Li19}.
	
	Let $G^{\rm ell}(x,y)$ denote the elliptic Green's function associated
	with the convex polygonal domain $\Omega$ (see also \eqref{G-cell}). It was proved in
	\cite[Lemma~3.2]{Ley09} and \cite[Lemma~2.1]{DDO20} that there exists
	$\eta\in(0,1)$, depending only on $\Omega$, such that
	\begin{align}
		\frac{|\nabla_x G^{\rm ell}(x,y) - \nabla_x G^{\rm ell}(\bar x,y)|}{|x-\bar x|^\eta}
		\leq 
		C \Big(|x-y|^{-N+1-\eta} + |\bar x-y|^{-N+1-\eta}\Big)
		\notag
	\end{align}
	for all $x,\bar x,y\in\Omega$ with $x\neq\bar x$ and
	$x,\bar x\neq y$. Moreover, by
	\cite[Proposition~1, Eq.~(7)]{Fromm1993book} and
	\cite[Eq.~(2.8)]{Ley09},
	\begin{equation}\label{eq:elliptic-Green-gradient}
		|\nabla_xG^{\rm ell}(x,y)|
		\le C|x-y|^{1-N}, \quad |\nabla_x \nabla_y G^{\rm ell}(x,y)|
		\le C|x-y|^{-N}
		\qquad x,y\in\Omega,\quad x\ne y.
	\end{equation}
	
	For fixed $j$ and $x_0$, choose a cutoff function $\chi_j$ such that $\chi_j = 1$ on $Q_j(x_0)$, $\operatorname{supp}\chi_j\subset Q_j'(x_0)$, satisfying the standard scaling laws
	\begin{align}\label{eq:cutoff1}
		|\nabla\chi_j|\le Cd_j^{-1},
		\quad
		|\partial_t\chi_j|+|\nabla^2\chi_j|\le Cd_j^{-2} .
	\end{align}
	We further arrange $\chi_j$ so that its spatial sections satisfy the following separation condition
	\begin{align}\label{eq:cutoff2}
		\operatorname{dist}\Big(\operatorname{supp}_x\big(\nabla\chi_j(t,\cdot)\big), Q_j(t,x_0) \Big)\geq d_j/{16} ,
		\quad
		\operatorname{dist}\Big(\operatorname{supp}_x\big(\nabla\chi_j(t,\cdot)\big), \Omega\setminus Q_j'(t,x_0) \Big)\geq d_j/{16} ,
	\end{align}
	where $Q_j(t,x_0):=\{ x: (t,x)\in Q_j \}$ is the spatial section of $Q_j(x_0)$ at time $t$. All constants in
	\eqref{eq:cutoff1}--\eqref{eq:cutoff2} are independent of $j$, $h$,
	and $x_0$.
	
	Set $\widetilde G:=\chi_jG$. Since $\partial_tG=\Delta G$, the product
	rule gives
	\begin{align}\label{eq:Delta-G-tilde}
		\Delta\widetilde G
		&=\chi_j\partial_tG
		+2\nabla\chi_j\cdot\nabla G
		+G\Delta\chi_j
		\notag\\
		&=(\chi_j\partial_tG
		-G\Delta\chi_j)
		+\nabla\cdot(2G\nabla\chi_j)
		\notag\\
		&=: \widetilde G_1 + \nabla\cdot \widetilde G_2
		 .
	\end{align}
	The Gaussian estimate \eqref{GausEst1}, the cutoff bounds
	\eqref{eq:cutoff1}, and H\"older's inequality on the spatial sections
	imply
	\begin{align}
		&\quad d_j^{-N}
		\| \widetilde G_1\|_{L^\infty_tL^1_x(\mathcal Q)}
		+
		\| \widetilde G_1\|_{L^\infty_tL^\infty_x(\mathcal Q)}
		\notag\\
		&
		+
		d_j^{-N-1} \| \widetilde G_2\|_{L^\infty_tL^1_x(\mathcal Q)}
		+
		d_j^{-1} \| \widetilde G_2\|_{L^\infty_tL^\infty_x(\mathcal Q)}
		\le C d_j^{-N-2}. \notag
	\end{align}
	
	We now derive a local gradient estimate using the elliptic Green
	representation formula:
	\begin{align}
		\widetilde G(t,z,x_0)
		&=
		\int_\Omega
		\widetilde G(t,y,x_0) \delta_z(y)  \,\mathrm{d}y
		=
		-
		\int_\Omega
		\widetilde G(t,y,x_0) \Delta_y G^{\rm ell}(z,y)  \,\mathrm{d}y
		\notag\\
		&
		=
		\int_\Omega G^{\rm ell}(z,y)
		\bigl(-\Delta_y\widetilde G(t,y,x_0)\bigr)\,\mathrm{d}y,
		\label{G-cell}
	\end{align}
	For $z\in Q_j(t,x_0)$, integration by parts
	with respect to $y$ yields
	\begin{align}
		|\nabla \widetilde G(t,z,x_0)|
		&=
		\bigg|\int_\Omega
		\nabla_zG^{\rm ell}(z,y)
		\Delta_y \widetilde G(t,y,x_0) \d y \bigg|
		\notag\\
		&\leq
		\bigg|\int_\Omega
		\nabla_zG^{\rm ell}(z,y)
		 \widetilde G_1(t,y,x_0) \d y 
		-
		\int_\Omega
		\nabla_z \nabla_y G^{\rm ell}(z,y)
		\cdot \widetilde G_2(t,y,x_0) \d y \bigg|
		\notag\\
		&\leq
		C
		\bigg|\int_\Omega
		|z-y|^{1-N}
		|\widetilde G_1(t,y,x_0)| \d y \bigg|
		+
		C
		\bigg|\int_\Omega
		|z-y|^{-N}
		|\widetilde G_2(t,y,x_0)| \d y \bigg|
		\notag\\
		&\leq
		C\bigl\||z-y|^{1-N}\bigr\|_
		{L^{1}(B_{d_j}(z))}
		\|\widetilde G_1(t,\cdot,x_0)\|_{L^\infty(\Omega)}
		\notag\\
		&\quad+
		C \bigl\||z-y|^{1-N}\bigr\|_
		{L^{\infty}(\Omega\setminus B_{d_j}(z))}
		\|\widetilde G_1(t,\cdot,x_0)\|_
		{L^1(\Omega)}
		\notag\\
		&
		\quad+
		Cd_j^{-N}
		\|\widetilde G_2(t,\cdot,x_0)\|_{L^1(\Omega)}
		\notag\\
		&\leq
		Cd_j d_j^{-N-2}
		+
		Cd_j^{1-N} d_j^{-2}
		+
		Cd_j^{-N} d_j^{-1}
		\leq
		C d_j^{-N-1}
		. \notag
	\end{align}
	In the estimate of the term involving $\widetilde G_2$, we used the
	separation property \eqref{eq:cutoff2}. Indeed, for each fixed $t$,
	$\widetilde G_2(t,\cdot,x_0)$ is supported in
	$\operatorname{supp}_x(\nabla\chi_j(t,\cdot))$, and hence $|z-y|\geq d_j/{16}$
	for all $z\in Q_j(t,x_0)$ and
	$y\in\operatorname{supp}_x(\nabla\chi_j(t,\cdot))$.
	
	Since $\widetilde G=G$ on $Q_j(x_0)$, we obtain
	\begin{align}\label{eq:grad-inf-inf}
		\|\nabla G\|_{L^\infty_tL^\infty_x(Q_j(x_0))}
		\leq
		Cd_j^{-N-1},
	\end{align}
	which proves \eqref{GFest04} for $k=0$.

	We next proceed by induction. For $k\geq1$, let $G_k := \partial_t^k G$.
	Arguing as in \eqref{eq:Delta-G-tilde}, we obtain
	\begin{align}
		\Delta (\chi_j G_k)
		=
		(\chi_j G_{k+1} - G_k \Delta \chi_j)
		+
		\nabla\cdot (2G_k\nabla\chi_j) . \notag
	\end{align}
	Applying the preceding argument to $G_k$ and using the Gaussian
	estimates for its time derivatives gives
    \vspace{-1.5mm}
	\begin{align}
		{\| \nabla \partial_t^k G(\cdot,\cdot,x_0) \|_{L^{\infty}L^\infty(Q_j(x_0))}
			\leq Cd_j^{-N-1-2k} } \qquad k=1,2,... . \notag
	\end{align}
	This completes the proof of \eqref{GFest04}; see also
	\cite[Lemma~4.1]{Li2017}.
	
	We now turn to the local H\"older estimate \eqref{GFest03}. Taking
	norms in \eqref{eq:Delta-G-tilde} and again using H\"older's
	inequality on the spatial sections, we find
	\begin{align}
		&\quad d_j^{-N}
		\| \Delta\widetilde G \|_{L_t^\infty L_x^1(\mathcal Q)}
		+
		\| \Delta\widetilde G \|_{L_t^\infty L_x^\infty(\mathcal Q)}
		\notag\\
		&\leq
		\| \chi_j \|_{L_t^\infty L_x^\infty(Q_j'(x_0))}
		\| \partial_t  G \|_{L_t^\infty L_x^\infty(Q_j'(x_0))}
			\notag\\
		&\quad
		+
		2
		\| \nabla \chi_j \|_{L_t^\infty L_x^\infty(Q_j'(x_0))}
		\| \nabla  G \|_{L_t^\infty L_x^\infty(Q_j'(x_0))}
		\notag\\
		&\quad
		+
		\| \Delta \chi_j \|_{L_t^\infty L_x^\infty(Q_j'(x_0))}
		\| G \|_{L_t^\infty L_x^\infty(Q_j'(x_0))}
		\notag\\
		&\leq
		C d_j^{-N-2} 
		+
		C d_j^{-1} d_j^{-N-1}
		+
		C d_j^{-2} d_j^{-N}
		\leq
		C d_j^{-N-2}
		\notag .
	\end{align}
	Here we have used \eqref{eq:grad-inf-inf} on the enlarged dyadic
	region. More precisely,
	\begin{align}
		\| \nabla  G \|_{L_t^\infty L_x^\infty(Q_j'(x_0))}
		&\leq
		\| \nabla  G \|_{L_t^\infty L_x^\infty(Q_{j+1}(x_0))} 
		+
		\| \nabla  G \|_{L_t^\infty L_x^\infty(Q_j(x_0))}
		\notag\\
		&\quad
		+
		\| \nabla  G \|_{L_t^\infty L_x^\infty(Q_{j-1}(x_0))}
		\leq
		C d_j^{-N-1}
		. \notag
	\end{align}
	
	For $z,\bar z\in\Omega$ with $z\neq\bar z$, the H\"older estimate for
	the elliptic Green's function gives
	\begin{align}
		&\frac{|\nabla_z\widetilde G(t,z,x_0)
			-\nabla_z\widetilde G(t,\bar z,x_0)|}{|z-\bar z|^\eta}
		=
		\bigg| \int_\Omega \frac{\nabla_z G^{\rm ell}(z,y) - \nabla_z G^{\rm ell}(\bar z,y)}{|z-\bar z|^\eta} \Delta_y \widetilde G(t,y,x_0) \d y \bigg|
		\notag\\
		&\leq C
		 \int_\Omega \Big(|z-y|^{-N+1-\eta} + |\bar z-y|^{-N+1-\eta}\Big) | \Delta_y \widetilde G(t,y,x_0) | \d y 
		\notag\\
		&\leq C
		\int_{|z-y|\leq d_j/2} |z-y|^{-N+1-\eta} | \Delta_y \widetilde G(t,y,x_0) | \d y
		+
		C
		\int_{|\bar z-y|\leq d_j/2} |\bar z-y|^{-N+1-\eta} | \Delta_y \widetilde G(t,y,x_0) | \d y
		\notag\\
		&\quad+ C
		\int_{|z-y|> d_j/2} |z-y|^{-N+1-\eta} | \Delta_y \widetilde G(t,y,x_0) | \d y
		+ C
		\int_{|\bar z-y|> d_j/2} |\bar z-y|^{-N+1-\eta} | \Delta_y \widetilde G(t,y,x_0) | \d y 
		\notag\\
		&\leq
		C
		\Big(\big\| |z-y|^{-N+1-\eta} \big\|_{L_y^{1}(B_{d_j}(z))} + \big\| |\bar z-y|^{-N+1-\eta} \big\|_{L_y^{1}(B_{d_j}(\bar z))} \Big)
		\| \Delta \widetilde G(t,\cdot,x_0) \|_{ L^\infty(\Omega)}
		\notag\\
		&\quad+
		C
		\Big(\big\| |z-y|^{-N+1-\eta} \big\|_{L_y^{\infty}(\R^N\backslash B_{d_j}(z))} + \big\| |\bar z-y|^{-N+1-\eta} \big\|_{L_y^{\infty}(\R^N\backslash B_{d_j}(\bar z))} \Big)
		\| \Delta \widetilde G(t,\cdot,x_0) \|_{L^1(\Omega)}
		\notag\\
		&\leq
		C
		d_j^{1-\eta}
		d_j^{-N-2}
		+
		C
		d_j^{-N+1-\eta}
		d_j^{-2}
		\leq
		C d_j^{-N-1-\eta} . \notag
	\end{align} 
	The resulting bound is uniform in $t$ and $x_0$.
	
	Since $\widetilde G=G$ on $Q_j(x_0)$, the definition of the
	restriction norm therefore gives
	\begin{align}
		\|G(\cdot,\cdot,x_0)\|_{L^\infty C^{1,\eta}(Q_j(x_0))}
		&\le
		\| \widetilde G(\cdot,\cdot,x_0)\|_
		{L^\infty L^\infty(\mathcal Q)}
		+
		\| \nabla \widetilde G(\cdot,\cdot,x_0)\|_
		{L^\infty L^\infty(\mathcal Q)}
		\notag\\
		&\quad
		+
		|\widetilde G(\cdot,\cdot,x_0)|_
		{L^\infty C^{1,\eta}(\mathcal Q)}
		\notag\\
		&
		\le C d_j^{-N} + C d_j^{-N-1} + C d_j^{-N-1-\eta}
		\le C d_j^{-N-1-\eta}
		\notag ,
	\end{align}
	where we have applied \eqref{GausEst1}, \eqref{eq:cutoff1} and \eqref{eq:grad-inf-inf} to estimate
	\begin{align}
		&\| \widetilde G(\cdot,\cdot,x_0)\|_
		{L^\infty L^\infty(\mathcal Q)}
		\leq
		\| G(\cdot,\cdot,x_0)\|_
		{L^\infty L^\infty(Q_j')}
		\le C d_j^{-N} ,
		\notag
	\end{align}
	and
	\begin{align}
		\| \nabla \widetilde G(\cdot,\cdot,x_0)\|_
		{L^\infty L^\infty(\mathcal Q)}
		\leq
		C d_j^{-1}
		\| G(\cdot,\cdot,x_0)\|_
		{L^\infty L^\infty(Q_j')}
		+
		C
		\| \nabla G(\cdot,\cdot,x_0)\|_
		{L^\infty L^\infty(Q_j')}
		\le C d_j^{-N-1} .
		\notag
	\end{align}
	This
	proves \eqref{GFest03}.
	
	It remains to prove \eqref{GFest01}. Using the representation formula
	(cf. \cite[Eq.~(4.16)]{Li19}) and integrating by parts in $y$, we have
	\begin{align}
		\partial_t^k H(t,x,x_0)
		=
		-\int_{\tau_{x_0}^h} \partial_t^k G(t,y,x) \nabla_y \tilde \delta_{x_0} \d y
		=
		\int_{\tau_{x_0}^h} \nabla_y \partial_t^k G(t,y,x) \tilde \delta_{x_0} \d y .
	\end{align}
	With the local high-order gradient estimate \eqref{GFest04} in hand, estimate \eqref{GFest01} follows by repeating the proof of
	\cite[Lemma~4.1]{Li19}, with $G$ therein replaced by $\nabla G$; see
	also \cite[Lemma~4.1]{Li2017}.
\end{proof}}

According to \cite[Lemma 5.1]{Li19}, we have the following result. 
\begin{lemma}\label{LocEEst} 
	{\it 
		Suppose that $\phi\in L^2(0,T;H^1_0(\Omega))\cap H^1(0,T;L^2(\Omega))$ and 
		$\phi_h\in H^1(0,T; S_h^0({\Omega}))$ satisfy the equation 
		\begin{align}\label{eq:err_eq_1}
			(\partial_t(\phi-\phi_h),\chi_h)
			+(\nabla (\phi-\phi_h),\nabla \chi_h) =0, 
			\quad\forall\, \chi_h\in  S_h^0(\Omega), \,\, \mbox{a.e.}\,\, t\in(0,T), 
		\end{align} 
		with $\phi(0)=0$ in $\Omega_j''$ and $0\leq j \leq J_*$.
		Then we have
		\begin{align} 
			\hspace{-8pt}
			&\vertiii{\partial_t(\phi-\phi_h)}_{Q_j} 
			+ d_j^{-1}\vertiii{\phi-\phi_h}_{1,Q_j} \notag\\ 
			& \leq 
			C\epsilon^{-3}\big(I_j(\phi_{h}(0))+X_j(I_h\phi-\phi) 
			+d_j^{-2} \vertiii{\phi-\phi_h}_{Q_j'}\big) \nn\\
			&\quad 
			+ (C h^{1/2}d_j^{-1/2} + C\epsilon^{-1}hd_j^{-1} + \epsilon) \big(
			\vertiii{\partial_t(\phi-\phi_h)}_{Q_j'}
			+d_j^{-1}\vertiii{\phi-\phi_h}_{1,Q_j'}\big), 
			\label{LocEngErr} 
		\end{align} 
		where
		\begin{align*}
			&I_j(\phi_{h}(0))=\|\phi_{h}(0)\|_{1,\Omega_j'} + d_j^{-1}\|\phi_{h}(0)\|_{\Omega_j'} \, ,\\[5pt]
			&X_j(I_h\phi-\phi)=
			d_j\vertiii{\partial_t(I_h\phi-\phi)}_{1,Q_j'}
			+\vertiii{\partial_t(I_h\phi-\phi)}_{Q_j'} 
			\nn\\
			&\qquad\qquad\qquad\,\, 
			+d_j^{-1}\vertiii{I_h\phi-\phi}_{1,Q_j'}+
			d_j^{-2}\vertiii{I_h\phi-\phi}_{Q_j'}  \, . 
		\end{align*}
		Here $\epsilon\in (0,1)$ is an arbitrary constant,  the positive constant 
	$C$ is independent of $h$, $j$ and $C_*$ (defined in \eqref{h-condition}),  and the norms $\vertiii{\cdot}_{k,Q_j'}$ and $\|\cdot\|_{k,\Omega_j'}$ are defined in \eqref{Q-norm}. 
	} 
\end{lemma}

\section{Discrete stochastic maximal $L^p$-$L^q$ regularity for $p\in(2,\infty)$, $q\in[2,\infty)$ or $p=q=2$}
\label{sec-4}
\setcounter{equation}{0}
We first give a pointwise-in-time bound for the continuous semigroup.
\begin{lemma}\label{lem:pt_est1}
	The following estimate holds for all $p\in(1,\infty)$ and any given $t>0$:
	\begin{align*}
		\| \nabla e^{t\Delta} g \|_{L^p}
		\leq
		C t^{-1/2}
			\| g \|_{L^p} .
	\end{align*}
\end{lemma}
\begin{proof}
	{For $g\in W^{1,p}_0(\Omega)$,} it
    suffices to show the norm equivalence
	$C^{-1} \|\nabla g \|_{L^p} \le \|(-\Delta)^{1/2} g\|_{L^p} \leq C \|\nabla g \|_{L^p}$ for $p\in(1,\infty)$.
	By \cite[Theorem 7.5(b)]{JK95}, we know
	\begin{align}
		\| (-\Delta)^{1/2} g \|_{L^{p}}\leq {C} \| \nabla g \|_{L^{p}} \notag
	\end{align}
	 for any $p\in(1,\infty)$.
	For the reverse norm bound, let the H\"older conjugate of $p$ be $p^\prime$. By using duality argument and $L^{p^\prime}$ Hodge decomposition on $\Omega$ (cf. \cite[Eq. (10.101)]{Gia13}, \cite[Corollary 2.7]{DDO20} and \cite[Corollary 3.10 and Remark 3.11]{Dauge92}), we derive
	\begin{align}
		\| \nabla g \|_{L^p}
		&\leq
		C
		\sup_{\| \nabla \phi \|_{L^{p^\prime}}\leq 1}
		\bigg| \int_\Omega \nabla g \cdot \nabla \phi dx \bigg|
		\notag\\
		&\leq
		C
		\sup_{\| \nabla \phi \|_{L^{p^\prime}}\leq 1}
		\bigg| \int_\Omega g \Delta \phi  dx \bigg|
		\notag\\
		&\leq
		C
		\sup_{\| \nabla \phi \|_{L^{p^\prime}}\leq 1}
		\| (-\Delta)^{1/2} \phi \|_{L^{p^\prime}}
		\| (-\Delta)^{1/2} g \|_{L^p}
		\notag\\
		&\leq
		C
		\| (-\Delta)^{1/2} g \|_{L^p} . \notag
	\end{align}
	Therefore,
	\begin{align*}
		\| \nabla e^{t\Delta} g \|_{L^p}
		\leq
		C
		\| (-\Delta)^{1/2} e^{t\Delta} g \|_{L^p}
		\leq
		C t^{-1/2}
		\| g \|_{L^p} ,
	\end{align*}
	where we have used $\| (-\Delta)^{1/2} e^{t\Delta} \|_{\mathcal L(L^{p},L^{p})}\le Ct^{-1/2}$ (cf. \cite[Proposition 2.2.15]{Lunardi12} and \cite[Eq. (1.2)]{Li19}).
\end{proof}

{We now proceed to prove a discrete counterpart of Lemma \ref{lem:pt_est1}.}
{\begin{lemma}\label{lem:pt_est2}
	We have the following estimate for all $p\in(1,\infty)$ and any given $t>0$:
	\begin{align}
		\| \nabla e^{t\Delta_h} g_h \|_{L^p}
		\leq
		C t^{-1/2}
		\| g_h \|_{L^p} \quad g_h\in S_h^0(\Omega) . \notag
	\end{align}
\end{lemma}
\begin{proof}
	We first consider the range $p\in[2,\infty)$. Let $g_h\in S_h^0$ be arbitrary, and let $g$ solve
	\begin{align}
		-\Delta g = -\Delta_h g_h\mbox{ in }\Omega,\qquad g=0\mbox{ on }\partial\Omega . \notag
	\end{align}
	By the definition of the discrete Laplacian, $g_h$ is the Ritz projection of $g$, that is,
	\begin{align}
		g_h = R_h g . \notag
	\end{align}
	Let $p'$ denote the H\"older conjugate of $p$. Using the norm equivalence in Lemma~\ref{lem:pt_est1}, the $L^{p'}$-duality argument for the Ritz projection, and the $W^{2,p'}$-regularity of the dual Dirichlet problem, we obtain
	\begin{align}
		&\| \nabla g \|_{L^p}^2
		\leq
		C
		\| (-\Delta)^{1/2} g \|_{L^p}^2
		\leq
		C
		\| g \|_{L^p} \| \Delta g \|_{L^p}
		\notag\\
		&
		\leq
		C
		(\| g_h \|_{L^p} + \| (1-R_h) g \|_{L^p}) \| \Delta_h g_h \|_{L^p}
		\notag\\
		&
		\leq
		C
		\| g_h \|_{L^p} \| \Delta_h g_h \|_{L^p}
		+
		C
		h \| \nabla g \|_{L^p} \| \Delta_h g_h \|_{L^p}
		\notag\\
		&
		\leq
		C
		\| g_h \|_{L^p} \| \Delta_h g_h \|_{L^p}
		+
		\epsilon C
		\| \nabla g \|_{L^p}^2
		+
		\epsilon^{-1} C
		h^2 \| \Delta_h g_h \|_{L^p}^2 . \notag
	\end{align}
	By absorption and the inverse inequality, we get
	\begin{align}
		\| \nabla g \|_{L^p}^2
		\leq
		C
		\| g_h \|_{L^p} \| \Delta_h g_h \|_{L^p} \notag
		.
	\end{align}
	The $W^{1,p}$-stability of the Ritz projection therefore yields the discrete Gagliardo--Nirenberg interpolation inequality
	\begin{align}
		\| \nabla g_h \|_{L^p}
		&=
		\| \nabla R_h g \|_{L^p}
		\leq
		C
		\| \nabla g \|_{L^p}
		\leq
		C
		\| g_h \|_{L^p}^{1/2} \| \Delta_h g_h \|_{L^p}^{1/2} . \notag
	\end{align}
	Using the
	analyticity of the discrete semigroup (cf. \cite[Theorem~2.1]{Li19}), we obtain
	\begin{align}\label{eq:p2-inf}
		\| \nabla e^{t\Delta_h} g_h \|_{L^p}
		\leq 
		C
		\| e^{t\Delta_h} g_h\|_{L^p}^{1/2} \| \Delta_h e^{t\Delta_h} g_h \|_{L^p}^{1/2}
		\leq
		C t^{-1/2}
		\| g_h \|_{L^p}\qquad\forall p\in[2,\infty) . 
	\end{align}
	
	We next treat the remaining range of integrability exponent. Motivated by the classical resolvent regularization of fractional powers (cf. \cite[Proposition~4.5]{Komatsu1966} and \cite[Proposition~2.1.3]{Haase2006book}) and by the mesh-dependent regularization introduced in \cite[Eq.~(4.15)]{LZ2026}, we define the regularized square-root holomorphic functional-calculus multiplier
	$$
		B_h(z) = \frac{z^{1/2}}{1+h^2 z} .
	$$
	With the following operator identity
	\begin{align}
		(-\Delta_h)^{1/2} - B_h(-\Delta_h)
		=
		(-\Delta_h)^{1/2} \big(I - (I+h^2 (-\Delta_h))^{-1}\big)
		=
		h^2 (-\Delta_h) B_h(-\Delta_h) ,  \notag
	\end{align}
	we derive
	\begin{align}\label{eq:norm-equiv}
		\| (-\Delta_h)^{1/2} g_h \|_{L^p}
		=
		\big\| \big(1 + h^2(-\Delta_h)\big) B_h(-\Delta_h) g_h \big\|_{L^p}
		\leq
		C \| B_h(-\Delta_h) g_h \|_{L^p} ,
	\end{align}
	where we have applied the following estimate, using the inverse inequality,
	\begin{align}
		\| h^2(-\Delta_h) f_h \|_{L^p}
		\leq
		C
		h^2 h^{-2} \| f_h \|_{L^p}
		=
		C\| f_h \|_{L^p}, \qquad\forall f_h\in S_h^0. \notag
	\end{align}
	
	For the corresponding continuous operator, we have the boundedness
	\begin{align}\label{eq:norm-equiv1}
		\| P_h B_h(-\Delta) g_h \|_{L^p}
		=
		\big\| P_h \big(I+h^2(-\Delta)\big)^{-1} (-\Delta)^{1/2} g_h \big\|_{L^p}
		\leq
		C \| (-\Delta)^{1/2} g_h \|_{L^p}
		\leq
		C \| \nabla g_h \|_{L^p} .
	\end{align}
	Here we have used the $L^p$-stability of $P_h$ and the uniform $L^p$-boundedness of the resolvent $\bigl(I+h^2(-\Delta)\bigr)^{-1}$. The latter follows from its Laplace-transform representation and the $L^p$-contractivity of the heat semigroup:
	\begin{align}
		\big\| \big(I + h^2(-\Delta)\big)^{-1} f \big\|_{L^p}
		=
		\Big\| \int_0^\infty e^{-t} e^{h^2 t \Delta} f \d t \Big\|_{L^p}
		\leq 
		C \| f \|_{L^p} . \notag
	\end{align}
	
	To compare the continuous and discrete multipliers, let
	\begin{align}
		w = (z+\Delta)^{-1} g_h, \qquad w_h = (z+\Delta_h)^{-1} g_h . \notag
	\end{align}
	Using the second resolvent identity and the relation
	$\Delta_hR_h=P_h\Delta$, we obtain
	(cf. \cite[p.~1261]{BP2017})
	\begin{align}
		w_h - P_h w
		&=
		\big((z+\Delta_h)^{-1} P_h - P_h (z+\Delta)^{-1}\big) g_h
		\notag\\
		&=
		(z+\Delta_h)^{-1}\big( P_h (z+\Delta) - (z+\Delta_h) P_h \big) (z+\Delta)^{-1} g_h
		\notag\\
		&=
		(z+\Delta_h)^{-1}( P_h \Delta - \Delta_h P_h ) w
		\notag\\
		&=
		(z+\Delta_h)^{-1} (-\Delta_h)P_h( 1-R_h ) w
		. \notag
	\end{align}
	Consequently, combining the $L^p$-stability of $P_h$ and $(-\Delta_h) (z+\Delta_h)^{-1}$ (see \cite{Li19}) with a duality argument for $R_h$ in $L^p$, the resolvent estimate, and the norm equivalence from Lemma \ref{lem:pt_est1}, we obtain the desired result. Here, the duality argument relies on the $W^{1,p'}$-stability of $R_h$ and the $W^{2,p'}$-regularity of the dual Dirichlet problem.
	\begin{align}
		\| w_h - P_h w \|_{L^p}
		&\leq
		C \| w - R_h w \|_{L^p}
		\leq
		C h \| \nabla w \|_{L^p}
		=
		C h \| \nabla (z+\Delta)^{-1} g_h \|_{L^p}
		\notag\\
		&
		\leq
		C h \| (z+\Delta)^{-1} (-\Delta)^{1/2} g_h \|_{L^p}
		\leq
		C h|z|^{-1} \| \nabla g_h \|_{L^p} . \notag
	\end{align} 
    Fix $\theta\in(\pi/2,\pi)$ and introduce the contour
$\mathcal C_\theta
=\{re^{i\theta}:r>0\}\cup\{re^{-i\theta}:r>0\}\subset\mathbb C$.
The upper ray is oriented toward zero and the lower ray
away from zero, and $z^{1/2}$ denotes the principal branch. Along this contour, $|B_h(re^{i\theta})|\leq C_\theta B_h(r)$.
	The Dunford integral representation
	(cf. \cite[Chapters~2 and~3]{Haase2006book}), together with the
	preceding resolvent error estimate, now gives
	\begin{align}\label{eq:norm-equiv2}
		&\big\| B_h(-\Delta_h) g_h - P_h B_h(-\Delta) g_h \big\|_{L^p}
		=
		\Big\| \frac{1}{2\pi i} \int_{\mathcal C_{\theta}} \frac{z^{1/2}}{1+h^2 z} (w_h - P_h w) \d z \Big\|_{L^p}
		\notag\\
		&
		\leq
		C h \| \nabla g_h \|_{L^p} \int_0^\infty \frac{r^{-1/2}}{1+h^2 r} \d r
		=
		C \| \nabla g_h \|_{L^p} \int_0^\infty \frac{s^{-1/2}}{1+s} \d s
		\leq C \| \nabla g_h \|_{L^p}.
	\end{align}
	Combining \eqref{eq:norm-equiv}--\eqref{eq:norm-equiv2}, we obtain
	\begin{align}\label{eq:norm-equiv3}
		\| (-\Delta_h)^{1/2} g_h \|_{L^p}
		\leq
		C \| \nabla g_h \|_{L^p} \qquad\forall p\in[2,\infty) .
	\end{align}
	
	We derive the reverse norm bound for the conjugate exponents by duality. Arguing as in Lemma~\ref{lem:pt_est1}, we obtain
	\begin{align}\label{eq:norm-equiv4}
		\| \nabla g_h \|_{L^{p'}}
		&\leq
		C
		\sup_{\| \nabla \phi \|_{L^{p}}\leq 1}
		\bigg| \int_\Omega \nabla  g_h \cdot \nabla \phi dx \bigg|
		\leq
		C
		\sup_{\| \nabla \phi \|_{L^{p}}\leq 1}
		\bigg| \int_\Omega g_h  \Delta_h R_h \phi  dx \bigg|
		\notag\\
		&\leq
		C
		\sup_{\| \nabla \phi \|_{L^{p}}\leq 1}
		\| (-\Delta_h)^{1/2} R_h \phi \|_{L^{p}}
		\| (-\Delta_h)^{1/2} g_h \|_{L^{p'}}
		\notag\\
		&\leq
		C
		\| (-\Delta_h)^{1/2} g_h \|_{L^{p'}} \quad\forall p'\in(1,2] ,
	\end{align}
	where we have used \eqref{eq:norm-equiv3} and the
	$W^{1,p}$-stability of $R_h$. It follows from \eqref{eq:norm-equiv4} and the analyticity of the discrete semigroup that
	\begin{align}\label{eq:p1-2}
		\| \nabla e^{t\Delta_h} g_h \|_{L^{p'}}
		\leq
		C
		\| (-\Delta_h)^{1/2} e^{t\Delta_h} g_h \|_{L^{p'}}
		\leq
		C
		t^{-1/2}
		\| g_h \|_{L^{p'}}
		 \quad\forall p'\in(1,2] .
	\end{align}
	The proof now is complete in view of \eqref{eq:p2-inf} and \eqref{eq:p1-2}.
\end{proof}}

We define the following two errors associated with the Green’s function: 
\begin{align}\label{err-green}
M(t,\cdot,x_0) = \Gamma_h(t,\cdot,x_0) - \Gamma(t,\cdot,x_0)
\quad\mbox{and}\quad
F(t,\cdot,x_0)=\nabla_{x_0} \Gamma_h(t,\cdot,x_0) -H(t,\cdot,x_0) . 
\end{align}
The lemma below, concerning the error $M$, was proved in  \cite[Lemma 4.4]{Li19}.
\begin{lemma}\label{LemGm22}
	{\it
		For any $x_0\in\Omega$, the functions $\Gamma_h(t,x,x_0)$, $\Gamma(t,x,x_0)$ and 
		$M(t,x,x_0)$ satisfy 
		\begin{align}
			&
			\| \partial_t M(\cdot,\cdot ,x_0)\|_{L^1((0,+\infty)\times\Omega)}  
			+\|t \partial_{tt} M(\cdot,\cdot , x_0)\|_{L^1((0,+\infty)\times\Omega)} \leq C ,
			\label{L1Ft22}\\
			&
			\|\partial_t\Gamma_h(t,\cdot, x_0)\|_{L^1(\Omega)}\leq Ce^{-\lambda_0t} , \qquad\forall\, t\ge 1, \label{L1Gammatx22}
		\end{align} 
		where the constants $C$ and {$\lambda_0>0$} are independent of $h$. 
	}
\end{lemma}

{We next state the following core lemma. Its proof, based on the local energy estimates in Section~\ref{sec:local}, is given in Section~\ref{ProofLm03}.}
\begin{lemma}\label{LemGm2}
	{\it
		For almost every $x_0\in \Omega$, the function $F(t,x,x_0)$ satisfies
		\begin{align}
			&
			\| F(\cdot,\cdot ,x_0)\|_{L^2(0,+\infty; L^1(\Omega))}  
			+
			h \| F(\cdot,\cdot ,x_0)\|_{L^\infty(0,+\infty; L^1(\Omega))}  
			\notag\\
			&\quad
			+
			h^{1+N/2} \| F(\cdot,\cdot ,x_0)\|_{L^\infty(0,+\infty; L^2(\Omega))}  
			+\|t^{1/2} \partial_{t}F(\cdot,\cdot , x_0)\|_{L^1((0,+\infty)\times\Omega)}
			\leq C ,
			\label{L1Ft}
		\end{align} 
		where the {constant} $C$ is independent of $h$ {and $x_0$}. 
	}
\end{lemma}

Thanks to the above preliminaries, we are in a position to prove the main result.

\begin{proof}[Proof of Theorem \ref{thm:disc_max_reg2}]
We differentiate \eqref{eq:stoh_conv_disc} w.r.t. $x_0$ and get
\begin{align}
	\nabla_{x_0} u_h(t,x_0,\omega) &=
	\int_0^t \int_\Omega \nabla_{x_0} \Gamma_h(t-s,x,x_0) f_h(s,x,\omega) \d x \d W_s
	\notag\\
	&=
	\int_0^t \int_\Omega F(t-s,x,x_0) f_h(s,x,\omega) \d x \d W_s
	\notag\\
	&\quad+
	\int_0^t \int_\Omega H(t-s,x,x_0) f_h(s,x,\omega) \d x \d W_s
	\notag\\
	&=: \mathcal M_h f_h + \mathcal K_h f_h , \notag
\end{align}
where $\mathcal M_h$ and $\mathcal K_h$ are the linear operators whose stochastic convolution kernels are $F$ and $H$, respectively. 

By the representation formula \eqref{expr-Gamma1} and the stochastic Fubini theorem (cf. \cite[Theorem 4.33]{DZ2014book}), it follows that 
\begin{align}
	\mathcal K_h f_h(t,x_0,\omega)
	&=
	\int_0^t \int_\Omega H(t-s,x,x_0) f_h(s,x,\omega) \d x \d W_s
	\notag\\
	&=
	\int_0^t \int_\Omega \bigg(\int_\Omega G(t-s,y,x) (-\nabla_y \tilde\delta_{x_0}(y)) \d y \bigg) f_h(s,x,\omega) \d x \d W_s
	\notag\\
	&=
	\int_\Omega (- \nabla_y \tilde\delta_{x_0}(y))  \bigg(  \int_0^t \int_\Omega  G(t-s,x,y) f_h(s,x,\omega) \d x \d W_s \bigg) \d y 
	\notag\\
	&=
	\int_\Omega \tilde\delta_{x_0}(y) \nabla_y \bigg(  \int_0^t \int_\Omega  G(t-s,x,y) f_h(s,x,\omega) \d x \d W_s \bigg) \d y .
\end{align}
For every fixed $t\in [0,\infty)$ and a.s. $\omega\in\O$, {according to \eqref{reg-Delta-est-2}, we can apply Schur's lemma \cite[p. 457]{Grafakos08} with the kernel therein replaced with $\tilde\delta_{x_0} (y)$ and get the boundedness}
\begin{align}\label{eq:K_h1}
	\| \mathcal K_h f_h(t,x_0,\omega) \|_{L^q_{x_0}}
	&\leq
	C
	\bigg\| \nabla_{x_0} \int_0^t \int_\Omega  G(t-s,x,x_0) f_h(s,x,\omega) \d x \d W_s \bigg\|_{L^q_{x_0}}, \quad q\in (1,\infty)
	.
\end{align}
The following mapping property of $\mathcal K_h$ results from \eqref{eq:CSMR}
\begin{align}\label{eq:K_h2}
		\| \mathcal K_h f_h(t,x_0,\omega) \|_{L^p_\omega L^p_t L^q_{x_0}}
		&\leq
		C_{p,q}
				\| f_h(t,x,\omega) \|_{L^p_\omega L^p_t L^{q}_{x}},
	\end{align}
	where $p\in(2,\infty)$, $q\in[2,\infty)$ or $p=q=2$, {and the time norm is taken over $[0,T]$}.

It remains to check the mapping property of $\mathcal M_h$. 
On the one hand, from the triangle inequality,
\begin{align}
	\| \mathcal M_h f_h (t,x_0,\omega) \|_{L^2_\omega L^2_t L^2_{x_0}}
	&\leq
	\| \mathcal K_h f_h (t,x_0,\omega) \|_{L^2_\omega L^2_t L^2_{x_0}}
	\notag\\
	&\quad+
	\bigg\| \int_0^t \int_\Omega \nabla_{x_0} \Gamma_h(t-s,x,x_0) f_h(s,x,\omega) \d x \d W_s \bigg\|_{L^2_\omega L^2_t L^2_{x_0}}.
\end{align}
Using the It\^o isometry, Fubini's theorem and the energy estimate, we derive
\begin{align}
	&\bigg\| \int_0^t \int_\Omega \nabla_{x_0} \Gamma_h(t-s,x,x_0) f_h(s,x,\omega) \d x \d W_s \bigg\|_{L^2_\omega L^2_t L^2_{x_0}}^2
	\notag\\
	&=
	\E \int_0^{{T}} \int_\Omega \int_0^t \bigg(\int_\Omega \nabla_{x_0} \Gamma_h(t-s,x,x_0) f_h(s,x,\omega) \d x \bigg)^2 \d s \d x_0 \d t
	\notag\\
	&=
	\E \int_0^{{T}} \int_{{s}}^{{T}} \int_\Omega \bigg(\nabla_{x_0} \int_\Omega \Gamma_h({t-s},x,x_0) f_h(s,x,\omega) \d x \bigg)^2 \d x_0 \d t \d s
	\notag\\
	&=
	\E \int_0^{{T}} \int_{{s}}^{{T}} \int_\Omega \bigg(\nabla_{x_0}\Big(e^{(t-s)\Delta_h} f_h(s,\cdot,\cdot)\Big)(x_0) \bigg)^2 \d x_0 \d t \d s
	\notag \\
	&\leq
	\E \int_0^{T} \int_\Omega f_h(s,x,\omega)^2 \d x \d s
	=
	\| f_h (t,x,\omega) \|_{L^2_\omega L^2_t L^2_{x}}^2 . \notag
\end{align}
We apply a similar argument to {\eqref{eq:K_h2} with $p=q=2$}
therein and get
\begin{align}
	\| \mathcal K_h f_h (t,x,\omega) \|_{L^2_\omega L^2_t L^2_{x}}
	\leq 
	C
	\| f_h (t,x,\omega) \|_{L^2_\omega L^2_t L^2_{x}} . \notag
\end{align}
Therefore $\mathcal M_h: L^2_\omega L^2_t L^2_{x} \rightarrow L^2_\omega L^2_t L^2_{x}$ is bounded, i.e.
\begin{align}\label{eq:Mh-bbd-L2}
	\| \mathcal M_h f_h (t,x,\omega) \|_{L^2_\omega L^2_t L^2_{x}}
	\leq 
	C
	\| f_h (t,x,\omega) \|_{L^2_\omega L^2_t L^2_{x}} .
\end{align}

On the other hand, by Lemma \ref{LemGm2}, it holds that
\begin{align}
\bigg\| \int_\Omega F(t-s,x,x_0)f(s,x,\omega)\d x \bigg\|_{L^\infty_\omega L^\infty_t L^\infty_{x_0} L^2_{s\leq t}}
&
\le 
\| F(t-s,x,x_0) \|_{ {L_t^\infty L_{x_0}^\infty} L^2_{s\leq t} L^1_x} \| f(s,x,\omega) \|_{L^\infty_\omega L^\infty_s L^\infty_x}
\notag\\
&\le C \| f(s,x,\omega) \|_{L^\infty_\omega L^\infty_s L^\infty_x}, \notag
\end{align}
where $L^p_{s\leq t}:=\{f(s), s\in [0,t]: f\in L^p([0,t])\}$, equipped with the $L^p([0,t])$ norm.
Second, {a similar argument to that used in proving \eqref{eq:Mh-bbd-L2} leads to}  
\begin{align}
\bigg\| \int_\Omega F(t-s,x,x_0)f(s,x,\omega) \d x \bigg\|_{L^2_\omega L^2_t L^2_{x_0} L^2_{s\leq t}}
\le  C\| f(s,x,\omega) \|_{L^2_\omega L^2_s L^2_x}. \notag
\end{align}

The interpolation of the above two inequalities implies that for any $p\in [2,\infty]$,
 $$
\bigg\| \int_\Omega F(t-s,x,x_0)f(s,x,\omega) \d x \bigg\|_{L^p_\omega L^p_t L^p_{x_0} {L^2_{s\leq t}}}
\le  C\| f(s,x,\omega) \|_{L^p_\omega L^p_s L^p_x}.
 $$
This, combined with the BDG inequality, yields that
\begin{align}
	\| \mathcal M_h f_h (t,x_0,\omega) \|_{L^p_\omega L^p_t L^p_{x_0}}^p
	&=
	\int_\Omega \E \int_0^T \bigg| \int_0^t \int_\Omega F(t-s,x,x_0) f_h(s,x,\omega) \d x \d W_s \bigg|^p \d t \d x_0
	\notag\\
	&\leq
	C_p^p
	\int_\Omega \int_0^T  \E\bigg( \int_0^t  \bigg( \int_\Omega F(t-s,x,x_0) f_h(s,x,\omega) \d x \bigg)^2 \d s \bigg)^{p/2} \d t \d x_0
	\notag\\
	&\leq
	C_p^p  \| f_h(s,x,\omega) \|_{L_\omega^{p} L_s^{p} L_x^{p}}^p. \notag
\end{align}
Together with \eqref{eq:K_h2}, we conclude
\begin{align}\label{eq:LpLp}
	\|  \nabla_x u_h (t,x,\omega) \|_{L^p_\omega L^p_t L^{p}_x}
	\leq
	C_p 
	\| f_h(s,x,\omega) \|_{L_\omega^{p} L_s^{p} L_x^{p}} ,
\end{align}
for all $p\in[2,\infty)$.

{In order to apply Proposition \ref{prop:stochastic-extrapolation} to extend the above $L^p$-$L^p$ result to $L^p$-$L^q$,
	we verify the hypotheses in Proposition \ref{prop:stochastic-extrapolation} as follows.
	In our setting, we fix $q\in[2,\infty)$ and let
	\[
	X=X_h:=(S_h^0(\Omega),\|\cdot\|_{L^q}),\qquad
	Y=Y_h:=(S_h^0(\Omega),\|\cdot\|_{W^{1,q}}),\qquad A=A_h:=-\Delta_h.
	\]
	\begin{enumerate}[label=\textup{(\roman*)},nosep]
		\item
		The embedding $Y_h\hookrightarrow X_h$ with $v_h\mapsto v_h$ and
		$v_h\mapsto(v_h,\nabla v_h)$ into scalar and vector-valued
		$L^q$ spaces gives uniform UMD and type $2$ constants.
		The space $\gamma(\mathbb R,X_h)$ is isomorphic to $X_h$ (see \cite[Remark 2.2]{Van12} and \cite[Lemma 2.1 and Corollary 2.2]{Van08}).
		
		\item
		$A_h=-\Delta_h$ is sectorial of angle less than $\pi/2$ on $X_h$,
		with uniform bounds w.r.t. $h$; see \cite[Theorem~3.1]{Li22}.
		
		\item
		Lemma~\ref{lem:pt_est2} and the Poincar\'e inequality yield
		$\|e^{-tA_h}\|_{\mathcal L(X_h,Y_h)}
		\leq C_q t^{-1/2}$, uniformly in $h$.
	\end{enumerate}
	Setting $p=q$ in \eqref{eq:LpLp} and using the Poincar\'e inequality, we obtain the base estimate for $S_{A_h}$ with $p_0=q$. Applying Proposition~\ref{prop:stochastic-extrapolation} to the zero extension of $f_h$ outside $[0,T]$, we then obtain
	\[
	\|\nabla u_h\|_{L^p(\O\times(0,T);L^q)}
	\leq
	C_{p,q}
	\|f_h\|_{L^p(\O\times(0,T);L^q)}
	\]
	for all $p\in(2,\infty)$ and $q\in[2,\infty)$, with a constant
	independent of $h$ and $T$. The endpoint case $p=q=2$ follows directly
	from \eqref{eq:LpLp}. Finally, letting $T\rightarrow\infty$, the norm relation
	\eqref{eq:norm-equiv3} completes the proof of
	Theorem~\ref{thm:disc_max_reg2}.
}

\end{proof}

\subsection{Proof of Lemma \ref{LemGm2}}
\label{ProofLm03}

In the remainder of this section, we apply Lemma \ref{LocEEst} to estimate $$\| F \|_{L^2 ((0,+\infty); L^1(\Omega))}
+
 h \| F \|_{L^\infty(0,+\infty; L^1(\Omega))}
+
h^{1+N/2} \| F \|_{L^\infty ((0,+\infty); L^2(\Omega))}  
+
\| t^{1/2} \partial_t F \|_{ L^1((0,+\infty)\times \Omega)} . $$
The proof consists of two parts. The first part concerns estimates for $t \in (0, 1)$, and the second part concerns estimates for $t\ge 1$, which is a consequence of the smoothing property of parabolic equations. 

\medskip

\emph{Step I: Local energy estimates on ${\mathcal Q}=(0,1)\times\Omega$.}

{After fixing the lower bound $M$ in the choice $C_*\in[M,2M)$, we take $h_0\leq1/(8M)$ in Theorem~\ref{thm:disc_max_reg2}.}
First, we present estimates in the domain ${\mathcal Q}=(0,1)\times\Omega$
with the restriction $h<1/(4C_*) $; see \eqref{h-condition}. In this case, the basic global energy estimates for \eqref{RGFdef2} and \eqref{RGFdef3} give
\begin{align} 
&\|\partial_t H \|_{L^2({\mathcal Q})}
+\|\partial_t\nabla_{x_0}\Gamma_h\|_{L^2({\mathcal Q})}
\leq C(\|H(0)\|_{H^1} +\|\nabla_{x_0}\Gamma_h(0)\|_{H^1})
\leq Ch^{-2-N/2} ,  \label{EstGamma-Q1} \\
&\|H\|_{L^\infty L^2({\mathcal Q})} 
+\|\nabla_{x_0}\Gamma_h\|_{L^\infty L^2({\mathcal Q})} 
\leq C(\|H(0)\|_{L^2} +\|\nabla_{x_0}\Gamma_h(0)\|_{L^2})
\leq Ch^{-1-N/2} ,\label{EstGamma-Q2}\\
&\|\nabla H\|_{L^2({\mathcal Q})} +\|\nabla_x \nabla_{x_0}\Gamma_h\|_{L^2({\mathcal Q})}
\leq C(\|H(0)\|_{L^2} +\|\nabla_{x_0}\Gamma_h(0)\|_{L^2})
\leq Ch^{-1-N/2} , \label{EstGamma-Q3} \\
&\|\partial_{tt}H\|_{L^2({\mathcal Q})}
+\|\partial_{tt}\nabla_{x_0}\Gamma_h\|_{L^2({\mathcal Q})}
\leq C(\|\Delta H(0)\|_{H^1} +\|\Delta_{h,x}\nabla_{x_0}\Gamma_h(0)\|_{H^1})
\leq Ch^{-4-N/2} ,  \label{EstGamma-Q4}
\end{align} 

where we have used \eqref{reg-Delta-est}, {the inverse inequality (Property (P2))} and \eqref{Detal-pointwise} to estimate $H(0)$ and $\nabla_{x_0}\Gamma_h(0)$, respectively.

 Hence, from \eqref{EstGamma-Q2}, {\eqref{h-condition}} and {H\"older's inequality on a dyadic block}, we have 
 \begin{align}\label{L2Gamma-Q}
\vertiii{H}_{Q_*} +\vertiii{\nabla_{x_0}\Gamma_h}_{Q_*} 
\le Cd_*\|H\|_{L^\infty L^2(Q_*)}+Cd_*\|\nabla_{x_0}\Gamma_h\|_{L^\infty L^2(Q_*)}
\le Cd_* h^{-1-N/2} \le CC_*h^{-N/2} .
\end{align}

Since the volume of $Q_j$ is $Cd_j^{N+2}$, we can decompose 
$$ \| F \|_{L^2 L^1(\mathcal Q)}
+
{h \| F \|_{L^\infty L^1(\mathcal Q)}}
+
h^{1+N/2} \| F \|_{L^\infty L^2(\mathcal Q)}  
+
\| t^{1/2} \partial_t F \|_{L^1 L^1(\mathcal Q)}
$$ 
in the following way:
\begin{align}\label{Bd31K2}
&\| F \|_{L^2 L^1(\mathcal Q)}
+
h \| F \|_{L^\infty L^1(\mathcal Q)}
+
h^{1+N/2} \| F \|_{L^\infty L^2(\mathcal Q)}  
+
\| t^{1/2} \partial_t F \|_{L^1 L^1(\mathcal Q)}  
\nn\\
&\leq  
\| F\|_{L^2 L^1(Q_*)}
+  
C h {d_*}^{N/2} \| F \|_{L^{\infty} L^2(Q_*)} 
+  
h^{1+N/2} \| F \|_{L^{\infty} L^2(Q_*)} 
+ d_* \| \partial_t F \|_{L^1 L^1(Q_*)} 
\notag\\
&\quad 
+
\sum_{j}\big(\| F \|_{L^2 L^1(Q_j)}
+
C {hd_j^{N/2} \| F \|_{L^{\infty} L^2(Q_j)}}
+  
h^{1+N/2} \| F \|_{L^{\infty} L^2(Q_j)} 
+ d_j \| \partial_t F \|_{L^1 L^1(Q_j)}  
\big)  \nn\\
&\leq Cd_{*}^{1+N/2} \big(d_*^{-1} \vertiii{F}_{Q_*} 
+ d_* \vertiii{\partial_t F}_{Q_*} 
\big) 
+\sum_{j}Cd_{j}^{1+N/2}
\big(d_j^{-1}\vertiii{F}_{Q_j} 
+ d_j \vertiii{\partial_t F}_{Q_j} 
\big) \nn\\
&\leq CC_*^{2+N/2}
+ {C}{\mathscr K} .
\end{align} 
{Here, in the second inequality, we have used the Poincar\'e inequality in the temporal direction for functions with zero average, the triangle inequality and H\"older's inequality to obtain
	\begin{align}
		\| F \|_{L^\infty L^2(Q_j)}
		&\leq 
		C
		d_j^{-2}
		\| F \|_{L^1 L^2(Q_j)}
		+
		C
		\| \partial_t F \|_{L^1 L^2(Q_j)}
		\notag\\
		&
		\leq 
		C
		d_j^{-1}
		\| F \|_{L^2 L^2(Q_j)}
		+
		C
		d_j
		\| \partial_t F \|_{L^2 L^2(Q_j)}
		, \notag
	\end{align}
	and, similarly,
	\begin{align}
		\| F \|_{L^\infty L^2(Q_*)}
		\leq 
		C
		h^{-1}
		\| F \|_{L^2 L^2(Q_*)}
		+
		C
		h
		\| \partial_t F \|_{L^2 L^2(Q_*)}
		. \notag
	\end{align}
}
In the last inequality, we have used  \eqref{EstGamma-Q1} and \eqref{L2Gamma-Q} to estimate
$$
Cd_{*}^{1+N/2} \big(d_*^{-1} \vertiii{F}_{Q_*} 
+ d_* \vertiii{\partial_t F}_{Q_*}
\big), 
$$ 
and introduced the notation 
\begin{align}\label{express-K}
	{\mathscr K}  :&=
	\sum_{j} d_{j}^{1+N/2}
	(
	d_j^{-1} \vertiii{F}_{Q_j}
	+
	\vertiii{F}_{1,Q_j}
	+ d_j \vertiii{\partial_tF}_{Q_j}
	)  . 
\end{align}

It remains to estimate ${\mathscr K} $. To this end,
we set ``$\phi_h=\nabla_{x_0}\Gamma_h$, $\phi=H$, $\phi_h(0)= - P_h \nabla_x \tilde\delta_{x_0}$ and $\phi(0)=-\nabla_x \tilde\delta_{x_0}$'' in Lemma \ref{LocEEst}.
Then we obtain 
\begin{align}\label{F1tE}
d_j^{-1} \vertiii{F}_{1,Q_j}+\vertiii{\partial_tF}_{Q_j} 
&\le 
C\epsilon^{-3}(\widehat{I_j}+\widehat{X_j} + d^{-2}_j \vertiii{ F }_{Q'_j} )  \\
&\quad 
+ (C h^{1/2}d_j^{-1/2} + C\epsilon^{-1}hd_j^{-1}  + \epsilon) (d_j^{-1}\vertiii{F }_{1,Q_j'} + \vertiii{\partial_{t}F }_{Q_j'}) \qquad \quad\,\, \nn
\end{align}
By using the interpolation error estimate (i.e. Property (P3)), the exponential decay estimates \eqref{Detal-pointwise}--\eqref{Detal-pointwise1} and Lemma \ref{GFEst1} (local energy estimates of regularized Green's function), we have

\vspace{-2mm}
\begin{align}
	\widehat{I_j}
	&:=
	\|P_h \nabla_x \tilde \delta_{x_0}\|_{1,\Omega_j'}
	+d_j^{-1} \|P_h \nabla_x \tilde \delta_{x_0}\|_{\Omega_j'} \notag\\
	&\le C \Big(\int_{\Omega_j'} h^{-2N-4}e^{-\frac {2|x-x_0|}{Kh}} \d x \Big)^{1/2}
	+
	C
	d_j^{-1} \Big( \int_{\Omega_j'} h^{-2N-2}e^{-\frac {2|x-x_0|}{Kh}}  \d x \Big)^{1/2}  \notag\\
	&\le C d_{j}^{N/2} ( h^{-N-2} +d_j^{-1} h^{-N-1}  ) e^{-c d_j/h} \notag\\
	&\leq C h d_{j}^{-3-N/2} 
	\\[5pt]
\widehat{X_j}
&:=
d_j \vertiii{\partial_t(I_h H-H)}_{1,Q_j'}
+\vertiii{\partial_t(I_h H - H)}_{Q_j'}  \nn\\
&\quad\, 
+ d_j^{-1} \vertiii{I_h H - H}_{1,Q_j'}
+ d_j^{-2} \vertiii{I_h H - H}_{Q_j'}  \nn\\
&
\leq Cd_j h \|\partial_{t} H\|_{L^2H^{2}(Q_j'')}
+ Cd_j^{-1} h \|H\|_{L^2H^{2}(Q_j'')} \nn\\
&
\leq C h d_j^{-3-N/2} \qquad\mbox{(\eqref{GFest01} is used with $\gamma=1$)} . \label{ovXj}
\end{align}

By substituting \eqref{F1tE}--\eqref{ovXj} into the expression of ${\mathscr K}$ in \eqref{express-K}, we have
\begin{align} 
{\mathscr K} 
&=\sum_{j} d_{j}^{1+N/2} 
(
d_j^{-1} \vertiii{F}_{Q_j}
+
\vertiii{F}_{1,Q_j}
+d_j \vertiii{\partial_tF}_{Q_j}
) 
\nn\\
&\leq 
\sum_{j} d_{j}^{N/2} \vertiii{F}_{Q_j}
+
C \sum_{j}  d_{j}^{1+N/2}
\epsilon^{-3}
\big({ h d_j^{-2-N/2}}
+ d_j^{-1} \vertiii{F}_{Q'_j} \big) \nn\\
&\quad
+ \sum_{j} (C h^{1/2}d_j^{-1/2} + C\epsilon^{-1}hd_j^{-1}  + \epsilon) d_j^{1+N/2}(\vertiii{F }_{1,Q_j'} + d_j \vertiii{\partial_{t}F }_{Q_j'}) \nn\\
&\leq 
 C_\epsilon C_*^{-1} + C_\epsilon  \sum_{j}  d_{j}^{N/2} \vertiii{F}_{Q_j'}  \nn
\\
&\quad 
+ \sum_{j} (C h^{1/2}d_j^{-1/2} + C\epsilon^{-1}hd_j^{-1}  + \epsilon) d_j^{1+N/2} (\vertiii{F }_{1,Q_j'} + d_j \vertiii{\partial_{t}F }_{Q_j'}),
\end{align} 
where $C_{\epsilon}$ is a positive constant depending on $\frac 1\epsilon$.

Since $\vertiii{F}_{Q_j'} \le C(\vertiii{F}_{Q_{j-1}} + \vertiii{F}_{Q_{j}} + \vertiii{F}_{Q_{j+1}})$, 
we can convert the $Q_j'$-norm in the inequality above to the $Q_j$-norm:
\begin{align}\label{slkll}
{\mathscr K} &\leq {C_{\epsilon}C_*^{-1}} 
+ C_\epsilon d_{*}^{N/2}  \vertiii{F}_{Q_*} 
+ C_\epsilon \sum_{j}  d_{j}^{N/2} \vertiii{F}_{Q_j} \nn\\ 
&\quad 
+
C( h^{1/2}d_{*}^{-1/2} + \epsilon^{-1}hd_{*}^{-1} + \epsilon)
 d_*^{1+N/2}
(d_*^{-1}\vertiii{F}_{Q_*} +\vertiii{F}_{1,Q_*}
+d_{*}\vertiii{\partial_tF}_{Q_*})
 \nn\\ 
&\quad 
+\sum_{j} (C h^{1/2}d_j^{-1/2} + C\epsilon^{-1}hd_j^{-1} + \epsilon) d_{j}^{1+N/2}
(d_j^{-1}\vertiii{F}_{Q_j} +\vertiii{F}_{1,Q_j}
+d_{j}\vertiii{\partial_tF}_{Q_j}
) \nn\\
&\leq 
{C_{\epsilon}C_*^{-1}} 
+C_\epsilon C_*^{2+N/2}+ 
C_\epsilon  \sum_{j} d_{j}^{N/2}  \vertiii{F}_{Q_j} 
+ (C C_*^{-1/2} + C_\epsilon C_*^{-1} + \epsilon) {\mathscr K} ,
\end{align}
where we have used $d_j\ge C_*h$, 
	\eqref{EstGamma-Q1}--\eqref{EstGamma-Q3} and \eqref{L2Gamma-Q} to estimate 
$$
\mbox{$\vertiii{F}_{Q_*}$,\,\, $\vertiii{F}_{1,Q_*}$ \,\, and \,\, 
$\vertiii{\partial_{t}F}_{Q_*}$}
$$
and used the expression of ${\mathscr K}$ in \eqref{express-K} to bound the terms involving $Q_j$.  
{By choosing $\epsilon$ small enough and then requiring the lower bound $M$ for $C_*\in[M,2M)$ to be sufficiently large (with its final choice made in Step III),} the last term on the right-hand side of \eqref{slkll} can be absorbed by the left-hand side. Therefore, we obtain
\begin{align} \label{K-L2F}
{\mathscr K}
&\leq
{C_{\epsilon}C_*^{-1}} +{C_{\epsilon}} C_*^{2+N/2}+ {C_{\epsilon}}
\sum_{j} d_{j}^{N/2} \vertiii{F}_{Q_j} .
\end{align}

\emph{Step II: The duality argument.} 
It remains to estimate $\vertiii{F}_{Q_j}$. We apply the parabolic duality argument. Let $w$ be the solution of the following backward parabolic equation on the domain $\mathcal Q$
\begin{align}
	\left\{\begin{array}{ll}
			-\partial_t w - \Delta w = v, \notag\\
		w(1)=0 ,  \notag \; 
	\end{array}\right.
\end{align}
with {homogeneous Dirichlet boundary conditions},  $\operatorname{supp}(v) \subseteq Q_j$ and $\| v \|_{L^2(Q_j)} = 1$.
After testing the above equation with $F$, we get
\begin{align}\label{eq:dual_test}
	[F, v] = (F(0), w(0)) + [\partial_t F, w] + [\nabla F, \nabla w] ,
\end{align}
where
\begin{align}
	(F(0), w(0)) 
	&= -(P_h \nabla_x \tilde \delta_{x_0} - \nabla_x \tilde \delta_{x_0}, w(0)) \notag\\
	&= -(P_h \nabla_x \tilde \delta_{x_0} - \nabla_x \tilde \delta_{x_0}, w(0) - I_h w(0))
	\notag\\
	&=-(P_h \nabla_x \tilde \delta_{x_0}-\nabla_x \tilde \delta_{x_0}, w(0) - I_h w(0))_{\Omega_j''}
	-(P_h \nabla_x \tilde \delta_{x_0} - \nabla_x \tilde \delta_{x_0}, w(0) - I_h w(0))_{(\Omega_j'')^c}
	\notag\\
	&=: \mathcal I_1 + \mathcal I_2 .
\end{align}
Notice that $\nabla_x \tilde \delta_{x_0}$ vanishes on $\Omega_j''$.
From the exponential decay estimate \eqref{Detal-pointwise}, \eqref{reg-Delta-est}, Properties (P2) and (P3), we know
\begin{align}
	| \mathcal I_1 | 
	\leq C h d_j^{-1-N/2} . \notag
\end{align}

{
To handle $\mathcal I_2$, we let \(W_j\) be a set containing \((\Omega_j')^c\) but its distance to \(\Omega_j\) is larger than \(C^{-1}d_j\). Since 
\[
|x - y| + s^{1/2} \geq C_1^{-1}d_j \quad \text{for } x \in W_j \text{ and } (s, y) \in Q_j
\]
for some positive constant \(C_1\), it follows that \((s, x) \in \bigcup_{k \leq j+ \log_2 C_1} Q_k(y)\) for \((s, y) \in Q_j\). Now, if we denote \(\widetilde{G}(\cdot, \cdot, y)\) as any extension of \(G(\cdot, \cdot, y)\) from \(\bigcup_{k \leq j+ \log_2 C_1} Q_k(y)\) to \(\mathcal{Q}\), then for \(x \in W_j\) we have \vspace{-1.5mm}
\begin{align}
	w(0, x) = \int_0^1 \int_\Omega G(s, y, x)v(s, y) \, \d y\d s &= \iint_{Q_j} G(s, x, y)v(s, y) \, \d y \d s
	\notag\\
	&= 
	\iint_{Q_j} \widetilde{G}(s, x, y)v(s, y) \, \d y\d s,
	\notag
\end{align}
where we have used the symmetry $G(s,y,x)=G(s,x,y)$ and {the fact that $v$ vanishes outside $Q_j$.}
Hence, we have
\begin{align}\label{eq:w(0)}
	\|w(0, \cdot)\|_{C^{1,\eta}(W_j)} 
	&\leq C \int_0^1 \int_{\Omega} \| \widetilde G(s, \cdot, y)\|_{C^{1,\eta}(W_j)} |v(s, y)| \, \d y \, \d s
	\notag\\
	&
	 \leq C \sup_{y \in \Omega} \| \widetilde G(\cdot, \cdot, y)\|_{L^\infty C^{1,\eta}(\mathcal Q)} \|v\|_{L^1(Q_j)}
	 \notag\\
	 &\leq
	 C
	 \sup_{y\in \Omega}\|G(\cdot,\cdot,y)\|_{L^{\infty}C^{1,\eta}(\cup_{k\le j+\log_2 C_1} Q_k(y))}\|v\|_{L^1(Q_j)}
	 ,
\end{align}
where the last inequality follows from taking the infimum over all possible extensions $\widetilde G$,  in accordance with the intrinsic characterization of Besov spaces given in  \cite{Rychkov99}, and using the identification $C^{1,\eta} \cong B^{1+\eta}_{\infty,\infty}$.}
Then,
\begin{align*}
&| \mathcal I_2 |  = |(P_h \nabla_x \tilde \delta_{x_0} - \nabla_x \tilde \delta_{x_0}, w(0) - I_h w(0))_{(\Omega_j'')^c}|\\
&\le \|P_h \nabla_x \tilde \delta_{x_0} - \nabla_x \tilde \delta_{x_0}\|_{L^1((\Omega_j'')^c)} \|w(0) - I_h w(0)\|_{L^{\infty}((\Omega_j'')^c)}\\
&\le C h^{\eta}\| w(0) \|_{C^{1,\eta}((\Omega_j')^c)}
\quad\mbox{((P3) and real interpolation inequality are used)}
\notag\\
&\le C
h^{\eta}
 \sup_{y\in \Omega}\|G(\cdot,\cdot,y)\|_{L^{\infty}C^{1,\eta}(\cup_{k\le j+\log_2 C_1} Q_k(y))}\|v\|_{L^1(Q_j)}
\quad\mbox{(\eqref{eq:w(0)} is used)}
\\
&\le C 
h^{\eta}
\sum_{k=0}^{j+\log_2 C} d_k^{-1-N-\eta} \| v \|_{L^1(Q_j)}\quad\mbox{(\eqref{GFest03} is used)}
\notag\\
&\le C h^{\eta} d_j^{-N/2-\eta} \| v \|_{L^2(Q_j)},
\end{align*}
where the parameter $\eta\in(0,1)$ is given in \eqref{GFest03}.

{Denote the index set $\mathcal I_*:=\{ 0,1,...,J_*, * \}$ and define the index mapping $\iota(i)=i$ for $i\in \mathcal I_*\setminus \{*\}$ and $\iota(i)=J_*+1$ for $i= *$.}
From \eqref{eq:dual_test}, we know
\begin{align}\label{FL2Qj} 
	\vertiii{F}_{Q_j} 
	&= \sup_{\vertiii{v}_{Q_j}=1} \bigg|(F(0), w(0)) + [\partial_t F, w] + [\nabla F, \nabla w] \bigg| \notag\\
	&\leq C  
	h^{\eta} d_j^{-N/2-\eta}
	+ \sup_{\vertiii{v}_{Q_j}=1} \bigg|[\partial_t F, (1-I_h)w] + [\nabla F, \nabla (1-I_h)w] \bigg| \notag\\
	&\leq 
	C  h^{\eta} d_j^{-N/2-\eta}
	+ 
	C\sum_{\substack{i\in\mathcal I_*\\|\iota(i)-j|\le 2}} (h^{2} \vertiii{F_t}_{Q_i} + h \vertiii{F}_{1, Q_i} ) 
	\| w \|_{L^2 H^2(Q_i^\prime)}
	\notag\\
	&\quad
	+ 
	C\sum_{\substack{i\in\mathcal I_*\\|\iota(i)-j|\geq 3}} (h^{1+\eta} \vertiii{F_t}_{Q_i} + h^\eta \vertiii{F}_{1, Q_i} ) 
	d_i^{N/2}
	\| w \|_{L^2 C^{1,\eta}(Q_i^\prime)} \notag\\
	&\leq 
	C h^{\eta} d_j^{-N/2-\eta}
	+ 
	C\sum_{\substack{i\in\mathcal I_*\\|\iota(i)-j|\leq 2}} (h^{2} \vertiii{F_t}_{Q_i} + h \vertiii{F}_{1, Q_i} ) 
	\notag\\
	&\quad
	+C\sum_{\substack{i\in\mathcal I_*\\|\iota(i)-j|\geq 3}} (h^{1+\eta}\vertiii{F_t}_{Q_i} 
	+ h^\eta \vertiii{F}_{1,Q_i})
	d_i^{1-\eta}
	\bigg(\frac{\min (d_i ,d_j )}{\max (d_i ,d_j )}\bigg)^{\eta+N/2} ,
\end{align}
{{where, in the last inequality, $\| w\|_{L^2 C^{1,\eta}(Q_i')}$ in the off-diagonal blocks can be estimated similarly to \cite[Eqs. (5.26) and (5.31)]{Li19}. To be more precise, for well-separated indices $|\iota(i)-j|\ge 3$, by the local H\"older estimate in Lemma \ref{GFEst1}, we have
\[
\|G(s-t,\cdot,y)\|_{C^{1,\eta}( Q_i'(t))}
\le C \max(d_i,d_j)^{-N-1-\eta},
\]
and, consequently, by the Green's representation formula,
\[
\begin{aligned}
	\|w(t,\cdot)\|_{C^{1,\eta}(Q_i'(t))}
	&=
	\Big\|\int_t^1\int_\Omega G(s-t,x,y)v(s,y) \d y \d s \Big\|_{C^{1,\eta}(Q_i'(t))}
	\notag\\
	&\le
	C \max(d_i,d_j)^{-N-1-\eta}
	\int_{Q_j}|v(s,y)|\,dy\,ds \\
	&\le
	C \max(d_i,d_j)^{-N-1-\eta}
	|Q_j|^{1/2}\|v\|_{L^2(Q_j)} \\
	&\le
	C d_j^{1+N/2}\max(d_i,d_j)^{-N-1-\eta}.
\end{aligned}
\]
Using H\"older's inequality on $Q_i'$, we obtain
\begin{align}
	&\|w\|_{L^2C^{1,\eta}(Q_i')}
	\leq
	C d_i
	\|w\|_{L^\infty C^{1,\eta}(Q_i')}
	\le
	C d_i d_j^{1+N/2}\max(d_i,d_j)^{-N-1-\eta} 
	\notag\\
	&
	=
	C d_i^{1-\eta-N/2} \frac{d_i^{\eta+N/2} d_j^{1+N/2}}{\max(d_i,d_j)^{N+1+\eta}} 
	\le
	C d_i^{1-\eta-N/2}\bigg(\frac{\min (d_i ,d_j )}{\max (d_i ,d_j )}\bigg)^{\min\{ 1+N/2, \eta+N/2\}} \notag\\
	&=
	C d_i^{1-\eta-N/2}\bigg(\frac{\min (d_i ,d_j )}{\max (d_i ,d_j )}\bigg)^{\eta+N/2} . \notag
\end{align}
}}

We also note the following summation formula for geometric series (cf. \cite[Eq. (5.33)]{Li19}):
\begin{align}\label{lalnu}
\sum_j d_{j}^{N/2} 
\bigg(\frac{\min (d_i ,d_j )}{\max (d_i ,d_j )}\bigg)^{\eta+N/2}
\leq 
Cd_i^{N/2} ,
\end{align}
for any $i\in\mathcal I_*$.

\emph{Step III: Closing the stability loop---boundedness of $\mathscr K$.}

By substituting \eqref{FL2Qj} into \eqref{K-L2F},  we obtain 
\begin{align*}%
{\mathscr K} &\leq 
C_{\epsilon}C_*^{-1}
+
C_\epsilon C_*^{2+N/2}
+
C_\epsilon
\sum_{j} d_{j}^{N/2} \vertiii{F}_{Q_j} \nn\\
&\leq
C_{\epsilon}C_*^{-1}
+
C_\epsilon C_*^{2+N/2}
+
C_\epsilon \sum_j \left(\frac{h}{d_j}\right)^{\eta}
\\
&\quad
+
C_\epsilon
\sum_{j}  d_{j}^{N/2} \sum_{\substack{i\in\mathcal I_*\\|\iota(i)-j|\leq 2}}
(h^{2}\vertiii{F_t}_{Q_i} 
+ h \vertiii{F}_{1,Q_i})
\notag\\
&\quad
+
C_\epsilon
\sum_j  d_{j}^{N/2}
\sum_{\substack{i\in\mathcal I_*\\|\iota(i)-j|\geq 3}} (h^{1+\eta}\vertiii{F_t}_{Q_i} 
+ h^\eta \vertiii{F}_{1,Q_i})
d_i^{1-\eta}
\bigg(\frac{\min (d_i ,d_j )}{\max (d_i ,d_j )}\bigg)^{\eta+N/2}
\notag\\
&\leq
C_{\epsilon}C_*^{-1}
+
C_\epsilon C_*^{2+N/2}
+
C_\epsilon \sum_j \left(\frac{h}{d_j}\right)^{\eta}
\\
&\quad
+
C_\epsilon
\sum_{*,i}  d_{i}^{N/2}
(h^{2}\vertiii{F_t}_{Q_i} 
+ h \vertiii{F}_{1,Q_i})
\notag\\
&\quad
+
C_\epsilon
\sum_{*,i}
d_i^{1-\eta+N/2} (h^{1+\eta}\vertiii{F_t}_{Q_i} 
+ h^\eta \vertiii{F}_{1,Q_i}) ,
\end{align*}
where, in the last inequality, we have changed the summation order and used \eqref{lalnu}.
We continue the estimate for $\mathscr K$ as follows:
\begin{align*}%
{\mathscr K}  
&\leq
C_{\epsilon}C_*^{-1}
+
C_\epsilon C_*^{2+N/2}
+
C_\epsilon C_*^{-\eta}
\notag\\
&\quad
+
C_\epsilon
\sum_{*,i} d_i^{1+N/2}\left(h \vertiii{F_t}_{Q_i} 
\bigg(\frac{h}{d_i}\bigg)
+ \vertiii{F}_{1,Q_i}\bigg(\frac{h}{d_i}\bigg) \right) \\
&\quad
+
C_\epsilon
\sum_{*,i}d_i^{1+N/2}\left(h \vertiii{F_t}_{Q_i} 
\bigg(\frac{h}{d_i}\bigg)^\eta
+ \vertiii{F}_{1,Q_i}\bigg(\frac{h}{d_i}\bigg)^\eta \right) \\
&\leq 
C_{\epsilon}C_*^{-1}
+
C_\epsilon C_*^{2+N/2}
+
C_\epsilon C_*^{-\eta}
+
C_\epsilon
d_{*}^{1+N/2}
\left(d_*\vertiii{F_t}_{Q_*} + \vertiii{F}_{1,Q_*} \right) \\
&\quad
+
C_\epsilon
\sum_{i} d_i^{1+N/2}\left(h \vertiii{F_t}_{Q_i} 
+ \vertiii{F}_{1,Q_i}\right)\bigg(\frac{h}{d_i}\bigg) \\
&\quad
+
C_\epsilon
\sum_{i} d_i^{1+N/2}\left(h \vertiii{F_t}_{Q_i} 
+ \vertiii{F}_{1,Q_i}\right)\bigg(\frac{h}{d_i}\bigg)^\eta \\
&\leq 
C_{\epsilon}C_*^{-1}
+
C_\epsilon C_*^{2+N/2}
+
C_\epsilon C_*^{-\eta}
+
C_\epsilon C_*^{-1} {\mathscr K}
+
C_\epsilon C_*^{-\eta} {\mathscr K}   .
\end{align*}
{By fixing $M$ sufficiently large, uniformly for $C_*\in[M,2M)$} 
($C_*$ is determined now), the terms involving ${\mathscr K}$ will be absorbed by the left-hand side of the inequality above. In this case, the inequality above implies 
\begin{align}
{\mathscr K} \le {C_{\epsilon}}.
\end{align}  
Substituting the last inequality into \eqref{Bd31K2} yields 
\begin{align}\label{L1FtQ1}
\| F \|_{L^2 L^1(\mathcal Q)}
+
{h \| F \|_{L^\infty L^1(\mathcal Q)}  }
+
h^{1+N/2} \| F \|_{L^\infty L^2(\mathcal Q)}  
+
\| t^{1/2} \partial_t F \|_{L^1 L^1(\mathcal Q)}  
&\leq  {C_{\epsilon}}.
\end{align}

\emph{Step IV: Long-time stability.}

For $t\in[1,+\infty)$ and any given $g_h\in S_h^0({\Omega})$, we define $f_h = (-\Delta_h)^{-1} g_h\in S_h^0(\Omega)$ and $f\in H_0^1(\Omega)$ to be the solution of the elliptic equation
\begin{align}
	-\Delta f= g_h . \notag
\end{align}
By the definition we have the relation $f_h = R_h f$ where $R_h$ is the Ritz projection associated to the domain $\Omega$. 
For $p\in(1,\infty)$, let $p'$ denote its H\"older conjugate.
The $W^{1,p}$-stability of the Ritz projection $R_h$
(cf. \cite[Theorem 4.1]{Schatz95} and \cite{LQXY24,LL21}),
together with elliptic regularity lift, yields
\begin{align}\label{eq:disc_moll_1}
	\| (-\Delta_h)^{-1} g_h \|_{W^{1,p}} 
	=\| f_h \|_{W^{1,p}} 
	\leq C \| f \|_{W^{1, p}}
	\leq C \| g_h \|_{W^{-1, p}} \qquad p\in(1,\infty) .
\end{align}
The Sobolev embedding $W^{1,p^\prime} \hookrightarrow L^{(p^\prime)^*}$ furthermore gives
\begin{align}
	\| (-\Delta_h)^{-1} g_h \|_{W^{1,p}} 
	\leq C \| g_h \|_{W^{-1, p}}
	\leq C \| g_h \|_{L^{q}} 
	\qquad p\in(N/(N-1), \infty)
	,  \notag
\end{align}
where $q=((p^\prime)^*)^\prime$, i.e. $\frac{1}{q} = 1 - \frac{1}{(p^\prime)^*} = 1 - \frac{1}{p^\prime} + \frac{1}{N} = \frac{1}{p} + \frac{1}{N}$. Similarly, if $p=\infty$, using \eqref{eq:elliptic-Green-gradient}, we get
\begin{align}\label{eq:disc_moll_2}
	&\| (-\Delta_h)^{-1} g_h \|_{W^{1,\infty}} 
	=\| f_h \|_{W^{1,\infty}} 
	\leq C \| f \|_{W^{1, \infty}}
	\notag\\
	&\quad
	\leq C \sup_{x\in\Omega} \Big| \int_\Omega \big(|G^{\rm ell}(x,y)| + |\nabla_x G^{\rm ell}(x,y)|\big) |g_h(y)| \d y \Big|
	\notag\\
	&
	\leq C \sup_{x\in\Omega} \Big| \int_\Omega \big(1 + |x-y|^{1-N}\big) |g_h(y)| \d y \Big|
	\leq C_\beta \| g_h \|_{L^{N+\beta}}
\end{align}
for any ${\beta} > 0$.

By applying \eqref{eq:disc_moll_2} once and \eqref{eq:disc_moll_1} $(N-1)$ times, we get
\begin{align}%
	\| (-\Delta_h)^{-N-1} g_h \|_{W^{1,\infty}} 
	\leq C \| (-\Delta_h)^{-N} g_h \|_{L^{N+\beta}}
	\leq C \| (-\Delta_h)^{-1} g_h \|_{L^{\frac{1}{1/(N+\beta)+(N-1)/N}}} 
	\leq C \| g_h \|_{L^{1}}
	. \notag
\end{align}
{The last inequality follows from
\begin{align}
	\| (-\Delta_h)^{-1} g_h \|_{L^{\frac{1}{1/(N+\beta)+(N-1)/N}}} 
	\leq
	C \| g_h \|_{W^{-1,\frac{1}{1/(N+\beta)+(N-1)/N}}}
	\leq
	C \| g_h \|_{L^1} , \notag
\end{align}
where we have used the Sobolev embedding $W^{1,p^\prime} \hookrightarrow L^{\infty}$ for any $p<\frac{N}{N-1}$.}
Consequently for any $p\in(1,\infty)$, we apply the above discrete smoothing effect in $W^{1,\infty}$ norm to get
\begin{align}\label{eq:disc_moll_4}
	&\| (-\Delta_h)^{-N-1} P_h \nabla_{x} \tilde \delta_{x_0} \|_{W^{-1, p}}
	\notag\\
	&=
	\sup_{\| v \|_{W_0^{1, p^\prime}} = 1} \big|\big((-\Delta_h)^{-N-1} P_h \nabla_{x} \tilde \delta_{x_0}, v\big)\big|
	=
	\sup_{\| v \|_{W_0^{1, p^\prime}} = 1} \big| \big(\tilde \delta_{x_0}, \nabla_{x} P_h (-\Delta_h)^{-N-1} P_h v\big)\big|
	\notag\\
	&\leq
	\sup_{\| v \|_{W_0^{1, p^\prime}} = 1} \| \tilde \delta_{x_0} \|_{L^1} \| \nabla_{x} P_h (-\Delta_h)^{-N-1} P_h v \|_{L^\infty}
	\leq
	C \sup_{\| v \|_{W_0^{1, p^\prime}} = 1} \| (-\Delta_h)^{-N-1} P_h v \|_{W^{1,\infty}}
	\notag\\
	&\leq
	C \sup_{\| v \|_{W_0^{1, p^\prime}} = 1} \| P_h v \|_{L^{1}}
	\leq
	C \sup_{\| v \|_{W_0^{1, p^\prime}} = 1} \| v \|_{L^{1}}
	\leq
	C ,
\end{align}
where we have used the $L^1$ stability of $P_h$.

Denote by $e^{t\Delta_{h}}$ the analytic semigroup generated by the linear operator $\Delta_{h}: S_h^0({\Omega}) \rightarrow S_h^0({\Omega})$. Then,
\begin{align}
	&\| \partial_t \nabla_{x_0} \Gamma_h (t,\cdot, x_0) \|_{L^2}
	\notag\\
	&=
	\| \partial_t e^{t\Delta_{h}} P_h \nabla_{x} {\tilde\delta_{x_0}} \|_{L^2} \notag\\
	&=
	t^{-3-N} \| (t\Delta_{h})^{3+N} e^{t\Delta_{h}} (-\Delta_{h})^{-2-N} P_h \nabla_{x} {\tilde\delta_{x_0}} \|_{L^2} \notag\\
	&\leq C
	t^{-3-N}\| (-\Delta_{h})^{-2-N} P_h \nabla_{x} {\tilde\delta_{x_0}} \|_{L^2} \notag\\
	&\leq
	C t^{-3-N} \| (-\Delta_{h})^{-1-N} P_h \nabla_{x} {\tilde\delta_{x_0}} \|_{W^{-1,2}} \notag\\
	&\leq C t^{-3-N} . \notag
\end{align}
Here, the first inequality follows from the analyticity of the
discrete semigroup, while the second and third inequalities
follow from \eqref{eq:disc_moll_1} and
\eqref{eq:disc_moll_4}, respectively, with $p=2$.
Similarly,
\begin{align}
	\| \nabla_{x_0} \Gamma_h(t,\cdot, x_0) \|_{L^2}
	=
	\| e^{t\Delta_{h}} P_h \nabla_{x} {\tilde\delta_{x_0}} \|_{L^2} 
	\leq C t^{-2-N} . \notag
\end{align}
Therefore from the energy decay estimates {(cf. \cite[Eqs. (5.36)--(5.39)]{Li19})}, {there exists $\lambda_0>0$ such that for any $t\ge 1,$}
\begin{align}
	&\| \partial_t \nabla_{x_0} \Gamma_h(t,\cdot,x_0) \|_{L^2} 
	+
	\| \partial_t H(t,\cdot,x_0) \|_{L^2} \notag\\
	&\leq 
	C e^{-\lambda_0 (t - 1)} 
	(\| \partial_t \nabla_{x_0} \Gamma_h(1,\cdot,x_0) \|_{L^2} 
	+
	\| \partial_t H(1,\cdot,x_0) \|_{L^2}
	)
	\leq 
	C e^{-\lambda_0 (t - 1)} , \notag\\
	&\| \nabla_{x_0} \Gamma_h(t,\cdot,x_0) \|_{L^2}
	+
	\| H(t,\cdot,x_0) \|_{L^2}
	\notag\\
	&\leq 
	C e^{-\lambda_0 (t - 1)} 
	(\| \nabla_{x_0} \Gamma_h(1,\cdot,x_0) \|_{L^2} 
	+
	\| H(1,\cdot,x_0) \|_{L^2} 
	)
	\leq 
	C e^{-\lambda_0 (t - 1)} . \notag
\end{align}

In view of the exponential decay estimates above and \eqref{L1FtQ1} from Step III, 
the proof of Lemma \ref{LemGm2} is complete.

\medskip

\def\baselinestretch{0.95}
{\small
\bibliographystyle{abbrv}
\bibliography{WMP}
}

\end{document}